\documentclass[11pt,leqno]{article}
\usepackage{epsfig}
\usepackage{amssymb}
\usepackage{latexsym}
\usepackage{fullpage}
\title{Stability of travelling-wave solutions for reaction-diffusion-convection
systems}
\author{E.C.M.Crooks}
\date{}\def\epsfsize#1#2{\ifnum#1>\hsize\hsize\else#1\fi}

 \newtheorem{definition}{Definition}[section]
  \newtheorem{theorem}[definition]{Theorem}
  \newtheorem{corollary}[definition]{Corollary}
  \newtheorem{lemma}[definition]{Lemma}
  \newtheorem{proposition}[definition]{Proposition}
\newcommand{\RN}{{\mathbb R}^N}
\newcommand{\upos}{u^{+}}
\newcommand{\umin}{u^{-}}
\newcommand{\utpos}{u^{T, +}}
\newcommand{\utmin}{u^{T, -}}
\newcommand{\upm}{u^{\pm}}

\newcommand{\pmxt}{(\pm x, t)}
\newcommand{\xt}{(x,t)}
\newcommand{\ups}{\dsol}
\newcommand{\fop}{{\mathcal{F}}}
\newcommand{\pdo}{\psi_{0}}
\newcommand{\pddo}{\psi_{1}}
\newcommand{\vpdo}{v^{\pdo}}
\newcommand{\dsol}{{\mathfrak{h}}}
\newcommand{\dsold}{\dsol^{\delta}}
\newcommand{\eg}{{\mathfrak{e}}}

\newcommand{\xzp}{x_{0}(\phi)}
\newcommand{\xop}{x_{1}(\phi)}
\newcommand{\N}{{\mathbb N}}
\newcommand{\rs}{{\mathbb R}}
\newcommand{\sub}{{\mathbf{s}}_{\eta, x_{0}}}
\newcommand{\super}{{\mathbf{S}}_{\eta, x_{1}}}
\newcommand{\argxz}{(x-x_{0} + \eta \alpha_{0} e^{- \nu t})}

\newcommand{\ent}{e^{- \nu t}}
\newcommand{\buc}{{\mathfrak {C}}}
\newcommand{\tw}{{\bf (TW)}}

\newcommand{\cont}{\buc^{1}}
\newcommand{\nup}{\nu^{+}}
\newcommand{\num}{\nu^{-}}
\newcommand{\eep}{{\mathfrak{e}}^{+}}
\newcommand{\eem}{{\mathfrak{e}}^{-}}
\newcommand{\eepm}{{\mathfrak{e}^{\pm}}}
\newcommand{\pz}{p_{0}}
\newcommand{\qz}{q_{0}}
\newcommand{\etaz}{\eta_{0}}
\newcommand{\pg}{{\mathfrak p}}
\newcommand{\pgi}{\tilde{\pg}_{i}}
\begin{document}
\newcommand{\rrn}{(\rs, \RN)}
\newcommand{\fcont}{C^1 (\RN \times \RN, \RN)}
\newcommand{\proof}{\noindent {\em Proof.~~}}
\newcommand{\fixt}{(\cdot, t)}
\newcommand{\m}{E^{-}}
\newcommand{\p}{E^{+}}
\newcommand{\nop}{\mathcal{N}}
\newcommand{\mop}{\mathcal{M}}
\newcommand{\epm}{E^{\pm}}
\newcommand{\pstar}{\phi^{*}}
\newcommand{\cpm}{C^{\pm}}
\newcommand{\bpm}{B^{\pm}}
\newcommand{\C}{{\mathbb{C}}}
\newcommand{\pt}{\phi_{\tau}}
\newcommand{\pmin}{\phi_{-2 \eta_{2}}}
\newcommand{\ppos}{\phi_{2 \eta_{2}}}
\newcommand{\kb}{\kappa_{b}}
\newcommand{\ptb}{\psi_{\tau, b}}
\newcommand{\vptb}{v^{\ptb}}
\newcommand{\hep}{h_{\epsilon}}
\newcommand{\qed}{\hfill $\Box$ \\[2mm]}
\newcommand{\pop}{{\mathcal {P}}}
\newcommand{\rop}{{\mathcal {R}}}
\newcommand{\what}{\hat{w}}
\newcommand{\hop}{{\mathcal{H}}}
\newcommand{\wstar}{w^{*}}
\newcommand{\ball}{B_{\buc^{1}}( \rho_{0})}
\newcommand{\nbd}{(-\delta_{0}, \delta_{0})}
\newcommand{\up}{\hat{y}}
\newcommand{\sep}{\sigma_{\epsilon}}
\newcommand{\mept}{M_{\frac{\epsilon}{2}}}
\newcommand{\gept}{\gamma_{\frac{\epsilon}{2}}}
\newcommand{\dt}{~ ds}
\newcommand{\tk}{\tilde{K}}
\newcommand{\ballep}{B_{\buc^{1}}( \nuep)}
\newcommand{\ts}{\tilde{s}}
\newcommand{\nep}{\nu_{\epsilon}}
\newcommand{\rep}{\rho_{\epsilon}}
\newcommand{\lop}{{\mathcal {L}}}
\newcommand{\gep}{\gamma_{\epsilon}}
\newcommand{\bucc}{\widetilde{\buc}}
\newcommand{\real}{~\mbox{Real}~}
\newcommand{\nuep}{\nu_{\epsilon}}
\newcommand{\kep}{K_{\epsilon}}
\newcommand{\tep}{T_{\epsilon}}
\newcommand{\half}{{\scriptstyle \frac{1}{2}}}
\newcommand{\yz}{y_{0}}
\newcommand{\tyz}{\tau (\yz)}
\newcommand{\yphi}{y}
\newcommand{\that}{\hat{t}}
\newcommand{\dep}{\delta_{\epsilon}}
\bibliographystyle{plain}
\maketitle

\begin{abstract}
We are concerned with the asymptotic behaviour of classical solutions of 
systems of the form
\begin{equation}
\label{abs}
\left\{ \begin{array}{l} u_{t} = A u_{xx} + f(u, u_{x}), ~~
x \in \rs, t>0, u(x,t) \in \RN, \\
u(x,0) = \phi (x),
\end{array} \right.
\end{equation}
where $A$ is a positive-definite diagonal matrix and $f$ is a ``bistable''
nonlinearity satisfying conditions which guarantee the existence of
 a comparison principle for
(\ref{abs}).  Suppose that (\ref{abs}) has a
 travelling-front solution $w$ with velocity $c$, that connects two 
stable equilibria of $f$. (There are hypotheses on $f$ under which such
a front is known to exist \cite{ET}.)
We show that if $\phi$ is bounded, uniformly continuously
differentiable
and  such that $\| w(x) - \phi (x) \| $ is small when
$|x|$ is large, then there exists $\chi \in \rs$ such that
\begin{equation}
\label{abst}
\| u( \cdot , t) - w( \cdot + \chi - ct) \|_{BUC^{1}} \rightarrow 0
~~\mbox{as}~~ t \rightarrow \infty.
\end{equation}
Our approach extends an idea developed by Roquejoffre, Terman and Volpert 
in the convectionless case, where $f$ is independent of $u_{x}$.
First $\phi$ is assumed to be increasing in $x$, and (\ref{abst})
proved via a homotopy argument. Then we deduce the result for
 arbitrary $\phi$
 by showing that there is an increasing function in 
the $\omega-$limit set of $\phi$.
\end{abstract}

\section{Introduction}

This paper is concerned with the asymptotic behaviour of classical solutions of the system

\begin{equation}
u_{t} = Au_{xx}+ f(u, u_{x}), ~~x \in \rs , ~t>0, ~ u(x,t) \in \RN,
\label{origeq}
\end{equation}

\begin{equation}
u(x,0)= \phi (x), ~~ x \in \rs,
\label{data}
\end{equation}
under the following hypotheses: \\

\noindent
{\bf (a)} $A$ is a positive-definite diagonal $N \times N$ matrix, \\ \\
$f : \RN \times \RN \rightarrow \RN$ is a continuously-differentiable
 function such that \\
 \begin{description}
 \item[(f1)] $f_{i}(q,p)=\tilde{f}_{i}(q_{1}, \ldots, q_{N}, p_{i}) $ ~~ (the $i$-th component of $f$ does not depend on $p_{j}$ for $j \neq i)$,
 \item[(f2)] $ {\displaystyle \frac{ \partial f_{i}}{\partial q_{j}}(q,p) >0}, ~~ i \neq j, ~ i,j = 1, \ldots, N, ~~ (q,p) \in \RN \times \RN,$
 \item[(f3)] $f(\m,0)=f(\p,0)=0$, where $\m < \p$, $\epm
 \in \RN$ and all the eigenvalues of $d_{q}f[\epm,0]$ 
 lie in the open left-half complex plane ({\em bistability condition}),
 \item[(f4)] there exists $\gamma \in (1,2)$ and an increasing function $\mu : [0, \infty) \rightarrow [0,\infty)$ such that
for each $p,q \in \RN$,
$$ \| f(q,p) \| \leq \mu (\|q\|)(1 + \| p \|^{\gamma}) ~~~
\mbox{($\| \cdot \|$ denotes a norm on $\RN$)}$$
 \end{description}
and \\ 
\tw ~there exists a monotone travelling-wave solution $w(x-ct)$ of (\ref{origeq})  such that $w(x) \rightarrow E^{\pm}$ as $x \rightarrow \pm \infty$, 
 and $w'(x)>0$ is bounded independently of $x$. (In fact, these properties of
$w$ together with the above hypotheses on $f$ ensure that $w'(x) \rightarrow
0$ at an exponential rate as $|x| \rightarrow \infty$. See the remark following
the proof of Lemma \ref{lemdecay}.)\\[1mm]

\noindent Note that \cite{ET} proves the existence of a wave $w$ satisfying \tw $~$
under hypotheses similar, though not identical, to {\bf (a), (f1)-(f4)},
together with an assumption on the nonexistence of stable equilibria of $f$
between $(\m, 0)$ and $(\p, 0)$. Such equilibria could prevent the
 existence of a front connecting $\m$ to $\p$ - see \cite{FM}.
For the {\em scalar} bistable {\em equation}
 (\ref{origeq}), in the convectionless case
when $f \in \rs$ and
 is independent of $u_{x}$, convergence to a travelling-front
solution $w$ from initial data $\phi$ 
is comprehensively treated in \cite{FM}. Stability of fronts for bistable
convectionless {\em systems} is developed in \cite{Vol3} and \cite{RTV}.
Here we extend this work  to nonlinearities dependent on $u_{x}$.

Throughout,  $\eg = (1, \ldots, 1)$ and $d_{q}f[q,p]$,
 $d_{p}f[q,p]$ denote the partial Fr\'{e}chet derivatives
 of $f$ at $(q,p) \in \RN \times \RN$ with
 respect to the first and second arguments of $f$ respectively. If
$q^{\pm} \in \RN$, then $q^{-} < ( \leq) q^{+}$ if $q^{-}_{i} < (\leq)
q^{+}_{i}$ for each $i \in \{1, \ldots, N\}$; $[q^{-}, q^{+}]$ denotes the 
set of $q \in \RN$ such that $q^{-} \leq q \leq q^{+}$.  For
$\Upsilon$  a subset of a  real or complex
vector space,  $k \in
 {\mathbb{ N}} \cup \{ \infty \}$, ${\mathfrak {C}}^{k}(\rs , \Upsilon)
 =BUC^{k}(\rs , \Upsilon )$, the space of
 functions $g : \rs \rightarrow \Upsilon$ 
such that $g$ and the derivatives
 of $g$ of order less than or equal to $k$ are bounded and uniformly
 continuous on $\rs$. For brevity, we write $\buc^{k} =
 {\mathfrak {C}}^{k} \rrn$
and $\bucc^{k} = {\mathfrak {C}}^{k} ( \rs, \C^{N})$.

 Known results yield,
under hypotheses {\bf (a)}, {\bf (f1) - (f4)}, that 
 there exists $\epsilon >0$
such that    system (\ref{origeq} - \ref{data})
  with initial data $\phi \in \buc^{1} (\rs, [\m - \epsilon \eg
 , \p + \epsilon \eg ])$ has a
 unique classical solution $u^{\phi}$ that exists for all time and
 depends continuously in $\buc^{1}$ on the initial data $\phi$.
 See the Appendix for references. We will prove that  if $\phi \in \buc^{1}$
is such that $\| w(x) - \phi(x)\|$ is small when $|x|$ is large, then
 $u^{\phi}$ converges to a shift of the travelling wave $w$,
 in the sense that there exists 
$\chi \in \rs$, depending on $\phi$, such that 

\begin{equation}
\| u^{\phi} ( \cdot , t ) - w( \cdot + \chi - ct) \|_{\buc^{1}} \rightarrow 0
~~\mbox{as}~~ t \rightarrow \infty.
\label{origconv}
\end{equation}
Let $v(x,t)=u(x+ct, t)$, where $u$ is a solution of (\ref{origeq}). Then 

\begin{equation}
v_{t} = A v_{xx} + cv_{x} + f(v,v_{x}).
\label{frameeq}
\end{equation}
 Note that $w$ is a stationary solution of (\ref{frameeq})
and that 
$v(x,0) = u(x,0)$ for all $x \in \rs$. 
 We seek $\chi \in \rs$ such that

\begin{equation}
\| v^{\phi} ( \cdot, t) - w(\cdot + \chi) \|_{\buc^{1}} \rightarrow 0 ~~\mbox{as}~~
t \rightarrow \infty.
\label{frameconv}
\end{equation}

\noindent ($v^{\phi}$ will denote the unique classical solution of
 (\ref{frameeq})
with initial data  $\phi \in \buc^{1}(\rs, [\m - \epsilon \eg
 ,\p + \epsilon \eg ])$ throughout.)

To prove (\ref{frameconv}), it will first be shown, in Theorem
 \ref{localstabthm}, that $w$ is
 ``locally''stable
in $\cont$; that is, given initial data $\phi$ which is a sufficiently small
$\buc^{1}$-perturbation of $w$, the corresponding solution $v^{\phi}$ of
 (\ref{frameeq}) converges in
$\cont$ to a translate of $w$ as $t \rightarrow \infty$.
This is a consequence of the fact that the spectrum of the linearisation
of (\ref{frameeq}) about $w$ is in a sector in the open left-half plane,
except for a simple eigenvalue at zero caused by the  translation invariance
of $(\ref{frameeq})$.
  For 
$g \in \buc^{2}$ define

\begin{eqnarray}
\lop g (x) &=& Ag''(x) + \{ c+ d_{p}f[w(x),w'(x)] \} g'(x) +
 d_{q}f[w(x),w'(x)]g
(x) \nonumber \\
& = & Ag''(x) + C(x) g'(x) + B(x) g(x),
\label{linear}
\end{eqnarray}
say; $B,C : \rs \rightarrow M^{N \times N}$ are 
uniformly continuous $N \times N$-matrix-valued functions
of $x$. 
Consider $\lop$ as an operator acting in $\buc$, with domain $\buc^{2}$. We abuse notation slightly by also using the symbol $\lop$ for the
complexification of $\lop$ when appropriate. 
The spectrum of $\lop$ is analysed in section 2. Section 3 is devoted
 to proving local stability of $w$ in $\cont$, following a method in
 \cite{Henry}.

 The main convergence 
result, Theorem \ref{mainthm}, is proved in two steps. First, in section 4,
$\phi \in \buc^{1}$ is assumed to be increasing, and convergent to
$\epm$ at $\pm \infty$ respectively.
 Our approach derives from that of \cite{Vol3}.
 A function $\phi^{*}$ is constructed
 from $\phi$ and the wave $w$
so that the solution $v^{\phi^{*}}$ of (\ref{frameeq}) corresponding to 
initial data $\phi^{*}$ satisfies (\ref{frameconv}).
 The corresponding result for $v^{\phi}$ is then deduced using a homotopy
 argument.
 Section 5  concludes the paper
by showing that for more general initial data $\phi$, close to $w$ at infinity,
there is an increasing function in the $\omega$-limit set of $\phi$. This
last step is motivated by \cite{RTV}. 
Note that the main convergence theorem Theorem \ref{mainthm} implies uniqueness
of travelling-front solutions of (\ref{origeq}) within a certain class
 - see Corollary \ref{uniquecor} for details. 

In an Appendix, we state some useful known results for (\ref{frameeq}) - 
namely a comparison principle, local/global existence theorems and
{\em a priori} bounds. Some wave-dependent sub- and super-solutions, 
useful in the stability analysis of $w$, are  also given. 
This material will often be referred to in the body of the paper.\\[3mm]

\noindent {\bf Acknowledgements} I am very grateful to J.F. Toland for
 discussions on the work of this paper. I would also like to thank the
referee for valuable comments and suggestions.
This work was carried out under the
support of EPSRC Research Grant Number GR/K96342.

\section{Properties of $\lop$}

Let $Y, W$ be   complex Banach spaces and let
 $L(Y,W)$ denote the space of bounded
linear operators from $Y$ into $W$. A linear operator ${\mathcal{A}
: \mathcal{D}(\mathcal{A})} \subset Y \rightarrow Y$ is said to be 
{\em sectorial in $Y$} if it is a closed densely-defined operator such
that for some $\omega \in \rs, \theta \in ( \frac{\pi}{2} ,  \pi),
M>0$,
$$
\Sigma = \{ \lambda \in \C : \lambda \neq \omega, | \arg (\lambda - \omega) |
< \theta \} \subset \rho ( {\mathcal{A}}),~~\mbox{the resolvent set of}~~{\mathcal{A}}, $$
and
$$ \| (\lambda I - {\mathcal{A}})^{-1} \|_{L(Y,Y)}
 \leq \frac{M}{|\lambda - \omega|}
~~\mbox{for all}~~ \lambda \in \Sigma, $$
(see \cite[p 33]{Lunardi}). 
 If $\mathcal{A}$ is sectorial in $Y$, then $\mathcal{A}$ is
the infinitesimal generator of an analytic semigroup $e^{t {\mathcal {A}}}$
 in the Banach space $Y$.

\begin{lemma}
\label{lemsectorial}
The operator $\mathcal{L}$ $
 : \buc^{2} \subset \buc \rightarrow \buc$
defined in (\ref{linear}) is sectorial in $\buc$.
\end{lemma}

\proof In (\ref{linear}), the matrices $A$ and $C(\cdot)$ are diagonal.
It follows from the scalar-valued-equation analysis of \cite[p 81, Corollary
3.1.9 (ii)]{Lunardi} that the operator $\mathcal{T}$ $ : \buc^{2}
 \subset 
\buc \rightarrow \buc$ defined by ${\mathcal{T}}g = A g'' + C(\cdot)
g'$ is sectorial. 

Define ${\mathcal{S}} : \buc \rightarrow \buc$ by ${\mathcal{S}} g = B(\cdot) g$. Clearly ${\mathcal{S}} \in L( \buc,  \buc)$. So
\cite[p 64, Proposition 2.4.1]{Lunardi} yields that $\mathcal{L}$ $= \mathcal{T}
$ $
+ \mathcal{S}, \mathcal{L}$ $ : \buc^{2} \subset \buc \rightarrow 
\buc$ is sectorial. \qed

For ${\cal{A} : \cal{D}(\cal{A})} \subset Y \rightarrow Y$
and  $Y^{\sharp}$ be a Banach space with ${\mathcal{D}(\mathcal{A})}
 \subset 
Y^{\sharp}$ and $ Y^{\sharp} \hookrightarrow Y$, where $\hookrightarrow$ 
denotes continuous embedding,  let  the
{\em part} of $\mathcal{A}$ in $Y^{\sharp}$\cite[p 40]{Lunardi}  be
${\mathcal{A}}^{\sharp}$, where
$$ \mathcal{D}( {\mathcal{A}}^{\sharp}) = \mbox{$\{ g \in$}
 \mathcal{D}(\mathcal{A})
: \mathcal{A}\mbox{$g \in Y^{\sharp}$} \} \subset Y^{\sharp}, ~~\mbox{and}~~
 {\mathcal{A}}^{\sharp} g = {\mathcal{A}}g ~~\mbox{for each
$g \in {\mathcal {D}}({\mathcal {A}}^{\sharp})$}.$$

\begin{lemma}
\label{lempartsectorial}
The part of $\mathcal{L}$ in $\buc^{1}$ is sectorial in $\buc^{1}$.
\end{lemma}

\proof Define ${\mathcal{M}} : \buc^{2} \subset \buc
\rightarrow \buc$ by ${\mathcal {M}}g = Ag''.$ The proof of Lemma
\ref{lemsectorial} shows that both $\mathcal{M}$ and $\mathcal{L}$ are
 sectorial in $\buc$. Let $\mu_{0} \in \rs$ be such that if $\mu \in \C$ 
and $\real \mu \geq \mu_{0}$, then given $f \in \buc$, 
$({\mathcal {L}}-\mu I)g = f$ and $({\mathcal {M}}- \mu I)h = f$ are 
 solvable for $g$ and $h$ respectively. 
Then, keeping in mind  that functions in $\bucc$ are vector-valued,
an argument similar to that in the proof of
 \cite[p 92, Proposition 3.1.18]{Lunardi} yields the existence of $K>0$, 
independent of $\mu \in \C$ with $\real \mu \geq \mu_{0}$,
such that 
$$ \| \mu ( \mu I - {\mathcal {L}})^{-1}
 \|_{L(\bucc^{1} , \bucc^{1} )} < K ~~\mbox{if}~~
\real \mu \geq \mu_{0}. $$
The result follows from \cite[p 43, Proposition 2.1.11]{Lunardi}. \qed

We turn now to the spectral analysis of $\mathcal{L}$. Denote the spectrum
of $\lop$ by $\sigma (\lop)$ and the essential spectrum by
 $\sigma_{ess}(\lop)$. (Here, as in \cite{Henry},  the essential spectrum 
of $\mathcal{L}$ is the complement, in $\sigma (\lop)$, of the set of those
eigenvalues of finite (algebraic) multiplicity\footnotemark[1]
 which are isolated points of
 $\sigma (\lop)$.)
\footnotetext[1]{An eigenvalue $\lambda_0$ which is an isolated point of
the spectrum is said to have finite (algebraic) multiplicity if
$\pop \buc$ is finite-dimensional, where $\pop$ is the linear operator
defined by $\pop = \frac{1}{2 \pi i} \int_{\partial \Omega} (\xi I -
\lop)^{-1}  ~d \xi$, $\Omega$ being a ball in $\C$, centre
$\lambda_0$,
such that $\sigma (\lop) \cap \bar{\Omega} = \{ \lambda_0 \}$
 \cite[p 181]{Kato}.} 
 Of crucial importance 
is the following lemma
concerning the eigenvalues of the ``asymptotic form of $\lop$ at infinity''.
 It makes critical use of the 
bistability condition {\bf (f3)}. We define

\begin{equation}
\cpm  =  \lim_{x \rightarrow \pm \infty} C(x)
=cI + d_{p}f[E^{\pm}, 0]  ~~~\mbox{and}~~~
\bpm  =  \lim_{x \rightarrow \pm \infty} B(x) =
d_{q}f[E^{\pm}, 0]. \label{defbpm}
\end{equation}

\begin{lemma}
\label{lemeigen}
Suppose that there exist $\tau \in \rs, \lambda \in \C$ and $z \in {\mathbb{C}}^{N}$
such that 
\begin{equation}
(-\tau^{2}A + i\tau C^{+}+ B^{+})z = \lambda z.
\label{eigeneqn}
\end{equation}
Then $\real \lambda <0$. The same conclusion holds if $C^{+}, B^{+}$ are replaced by $C^{-}, B^{-}$ in (\ref{eigeneqn}).
\end{lemma}

\proof By condition {\bf (f3)}, all the eigenvalues of $\bpm$ lie
 in the open left-half complex plane. By condition {\bf (f1)}, $\cpm$ 
are diagonal and by condition {\bf (f2)}, $\bpm$ each have positive
off-diagonal elements. So the
result follows immediately from \cite[p 234, Lemma 4.1]{Vol3}. 
\qed

\begin{lemma}
\label{lemessential}
$\sigma_{ess}(\lop) \neq \emptyset$, and
 there exists $\beta >0$ such that if $\lambda  \in \sigma_{ess}(\lop)$
 then $\real \lambda < - \beta$.
\end{lemma}

\proof Let
\begin{equation}
\label{curves}
S^{\pm} = \{ \lambda \in \C : \mbox{det}~(- \tau^{2}A + i\tau \cpm + \bpm
-\lambda I) = 0~~\mbox{for some}~~ \tau \in \rs \}.
\end{equation}
Then Lemma \ref{lemeigen} shows that
$$ \lambda \in S^{+} \cup S^{-} \Rightarrow 
\real \lambda < 0.$$
\cite[p 140, Theorem A.2]{Henry} yields that $S^{\pm}$ each consists
 of a finite
number of algebraic curves parametrised by a real number $\sigma$, which are
asymptotically parabolic : $\lambda (\sigma) = - \sigma^{2} \alpha + O(\sigma)$
as $\sigma \rightarrow \infty$, where $\alpha$ is on the diagonal of $A$.
Moreover, $S^{+} \cup S^{-} \subset \sigma_{ess}(\lop)$ and $\sigma_{ess}(\lop)
\subset \Lambda$, where $\C \backslash \Lambda$ is the component of
$\C \backslash (S^{+} \cup S^{-})$ which contains the right-half plane.

Since $S^{\pm}$ are contained in the open left-half plane, $\Lambda$ is
also. Moreover, $S^{\pm}$ each consist of a finite number of algebraic curves
parametrised by $\sigma$, the real parts of which tend to $-\infty$ as 
$\sigma \rightarrow \pm \infty$. Whence $\Lambda$ is bounded away from
the imaginary axis. The result follows. \qed
 
We next show, using Lemma \ref{lemeigen}, 
 that the bistability condition {\bf (f3)} implies that bounded
solutions of certain equations must  decay at infinity.

\begin{lemma}
\label{lemdecay}
Suppose that there exist $\lambda \in \C, \real \lambda \geq 0$ and
$g \in \bucc^{2}$ such that $\lop g = \lambda g + \psi_{0}$, where
$\psi_{0} \in \buc$ is such that $\psi_{0} (x) \rightarrow 0$ as $|x|
\rightarrow \infty$. Then $\| g(x) \| \rightarrow 0$ as $|x| \rightarrow 
\infty$. If $\psi_{0} \equiv 0$, then there exist $M, \omega >0$ such that 
$\| g(x) \| \leq Me^{-\omega |x|}$ for all $x \in \rs$.
\end{lemma}

\proof  Define  $ \hat{h} = \left( \begin{array}{c} g \\ g' \end{array} 
\right)$, 
 $M^{+}=  \left(
\begin{array}{cc}
0 & I \\
-A^{-1}\{B^{+} - \lambda I \} & -A^{-1}C^{+} 
\end{array}
\right)$, $\hat{r} (x) = \hat{H}(x) \hat{h}(x) +
 \left( \begin{array}{c} 0 \\ \psi_{0} (x) \end{array} \right)$,
where   
$\hat{H} (x) =  \left(  
\begin{array}{cc}
0 & 0 \\
-A^{-1}\{ B(x) - B^{+}\} & -A^{-1}\{C(x)- C^{+}\} 
\end{array}
\right) $.
Then $ \hat{h}'(x)= M^{+}\hat{h}(x) + \hat{r}(x), ~x \in \rs$,
where $\hat{h}$ is bounded on $\rs$, and $\hat{r}(x) \rightarrow 0$
as $x \rightarrow \infty$. By Lemma \ref{lemeigen}, $M^{+}$ has no purely
imaginary eigenvalues. So, as in the proof of \cite[p 330, Theorem 4.1]{CL}, 
there exist $K, \alpha, \sigma > 0$, a real nonsingular matrix
$P \in M^{2N \times 2N}$ and operators $U_{1}(t), U_{2}(t)$ such that
\begin{equation}
\| U_{1}(x) \|  \leq  K e^{-\alpha x}, ~~ x \geq 0 ~~\mbox{and}~~
\label{decay1} 
\| U_{2}(x) \|  \leq  K e^{\sigma x}, ~~ x \leq 0,
\end{equation}
and $h = P \hat{h}$, $r = P \hat{r}$ satisfy
$$ h(x) = U_{1}(x) h(0) + U_{2}(x) k + 
\int_{0}^{x}U_{1}(x-s)r(s)~ds - \int_{x}^{\infty} U_{2}(x-s)r(s)~ds $$
where $$ k = h(0) + \int_{0}^{\infty}U_{2}(-s)r(s)~ds.$$
Estimates (\ref{decay1}) together with the facts that
$\hat{h}$ is bounded and $\hat{r} \rightarrow 0$ as $ x \rightarrow \infty$
 yield
that $h(x) \rightarrow 0$ as $x \rightarrow \infty$.
The exponential decay in the case $\psi_{0} \equiv 0$ follows from the
proof of \cite[p 330, Theorem 4.1]{CL}. \qed

\noindent {\bf Remark.} Clearly, due to translation invariance,
$ \lop w' =0, $
where $w'$ is the derivative of the travelling wave $w$. By hypothesis \tw,
$w'$ is bounded on $\rs$, so Lemma \ref{lemdecay} yields that $w'$ decays
exponentially to zero at $\pm \infty$. Also by \tw,  $w'(x) >0$
for all $x \in \rs$. Thus $\lop u =0$ has a positive solution which decays
exponentially to zero at infinity. Further, Lemma \ref{lemessential}
shows that zero is not in the essential spectrum of $\lop $, so it must be  
an isolated point of the spectrum and an eigenvalue of
 finite multiplicity. 

\begin{lemma}
\label{lemkrein}
\begin{itemize}
\item[(i)] For $\lambda \in \C \backslash \{0\}$ with $\real \lambda \geq 0$,
there are no non-zero solutions of the equation
\begin{equation}
\label{kreineigen}
\lop g = \lambda g, ~~ g \in \bucc^{2}.
\end{equation}
\item[(ii)] Let $ g \in \bucc^{2}$ be a solution of $\lop g =0$.
 Then there exists
$k \in \rs$ such that $g = kw'$.
\end{itemize}
\end{lemma}

\proof We aim to apply \cite[p 208, Theorem 5.1]{Vol3}. For this, note
that {\bf (f2)} and {\bf (f3)} imply that the matrix is irreducible in
the functional sense (defined in \cite[p 208]{Vol3}); this follows
from {\bf (f2)} alone when $N \geq 2$.
Now  let $\lambda \in \C, ~ \real \lambda \geq 0$, and suppose that
$g \in \bucc^{2}$ satisfies $\lop g = \lambda g$. That $\| g(x) \|
\rightarrow 0$ as $|x| \rightarrow \infty$ follows from Lemma
\ref{lemdecay} with $\psi_{0} \equiv 0$.
 The remark preceding this theorem together 
with \cite[p 208, parts (1) and (2) of Theorem 5.1]{Vol3} then yield
{\bf (i)} and {\bf (ii)}. \qed

\begin{proposition}
\label{propnegspec}
There exists $\gamma > 0$ such that if $\lambda \in \C$ belongs to
$\sigma ({\mathcal {L}}) \backslash \{ 0 \}$, then $\real \lambda < - \gamma$.
\end{proposition}

\proof Lemma \ref{lemessential} and Lemma \ref{lemkrein} {\bf (i)} show that
any non-zero point of $\sigma ({\mathcal {L}})$ lies in the open left-half complex plane.
If there is a sequence $\{ \lambda_{n} \} \subset \sigma ({\mathcal {L}})
 \backslash \{0\}$
such that $\real \lambda_{n} \uparrow 0$ as $n \rightarrow \infty$, then
by Lemma \ref{lemsectorial}, $\{ \mbox{Imag}~ \lambda_{n} \}$ is bounded.
Whence there is a subsequence $\{ \lambda_{k} \}$ and $\mu \in \sigma 
({\mathcal {L}})$ (a closed set), $\real \mu =0$, such that $\lambda_{k}
\rightarrow \mu$ as $k \rightarrow \infty$. But this contradicts Lemma
\ref{lemessential}. \qed

 Lemma \ref{lemkrein} {\bf (ii)}
 shows that the nullspace of $\lop$ is one-dimensional. We need additional
 information  to exploit this.
%
%
%
%
Recall that zero is an isolated eigenvalue of $\lop $. Let $\Omega$ denote a ball
in $\C$ with centre zero such that $\sigma (\lop ) \cap \overline{\Omega} = \{ 0 \}$.
Then for $\lambda \in \partial \Omega$, $(\lambda I - \lop )^{-1} : \buc 
\rightarrow \buc$ is a bounded linear operator; a bounded
linear operator $\pop$ is defined by 
\begin{equation}
\label{contour}
\pop = \frac{1}{2 \pi i} \int_{\partial \Omega} (\xi I - \lop ) ^{-1} ~d \xi,
\end{equation}
(see \cite[p 178]{Kato} or \cite[p 402]{Lunardi}). 
Let $X= \buc, X_{1} = \pop X$ and $X_{2}= (I-\pop )X$. \label{decomp}
 \cite[p 30, Theorem 1.5.2]{Henry} and \cite[p 402, Proposition A.1.2]{Lunardi} show that $\pop$
 is a projection, $X = X_{1} \oplus X_{2}$ and
 $\pop X$ is a subset of the domain of $\lop^{n}$ for 
each $n$. Moreover,   if $\lop_{j}$ is the restriction
of $\lop$ to $X_{j} \cap \buc^{2}$, then
\begin{eqnarray*}
& \lop_{1}&: X_{1} \rightarrow X_{1} ~~\mbox{ is bounded},
~~~~\sigma (\lop_{1}) = \{ 0 \} ~~\mbox{and} \\
&\lop_{2}& : X_{2} \cap \buc^{2} \subset X_{2} \rightarrow X_{2}, 
~~ \sigma (\lop_{2}) = \sigma (\lop )
 \backslash \{0 \} ~~\mbox{$(\neq \emptyset$, by
Lemma \ref{lemessential})}.
\end{eqnarray*}
 Note that since
$\pop , I-\pop $ are bounded operators by definition, $X_{1}$ and $X_{2}$ are
 closed subspaces of $X$.

\begin{lemma}
\label{lemprojection}
$X_{1} = \mbox{span} \{ w' \}$ and there exists $w^{*} \in X^{*}$ such that
\begin{equation}
\pop g = w^{*} (g) w' ~~\mbox{for each}~~ g \in X,\label{proj} ~~\mbox{and}~~ 
w^{ *} (w') = 1.
\label{projnormal}
\end{equation}
\end{lemma}

\proof \cite[p 405, Proposition A.2.2]{Lunardi} shows that ker $\lop \subset
 X_{1}$. 
Since $0 \not\in \sigma_{ess}(\lop)$,  $X_{1}$ is
finite-dimensional (see the footnote following the definition of  
$\sigma_{ess}(\lop)$). So $\sigma (\lop_{1})$ consists entirely of eigenvalues,
the number of which, counted according to algebraic multiplicity, equals the
dimension of $X_{1}$. 
It is shown in \cite[p 210,  proof of Theorem 5.1 (3)]{Vol3}
that Range $~\lop \cap ~\mbox{span}~\{ w' \}=0$.
 Thus zero is an eigenvalue of $\lop_{1}$
of multiplicity one, whence ker $\lop = X_{1}$. Since $\pop$ is a bounded projection,
the existence of $w^{*}$ as in the statement of the lemma follows. \qed
 
We will need two estimates on the behaviour of ${\mathcal{L}}_{2}$.
Define $ \gamma_{0} = - \sup~ \{ \real z : z \in \sigma (\lop_{2})\}$.
 By \mbox{Proposition
\ref{propnegspec}}, $\gamma_{0} >0$.

\begin{lemma}
\label{lemc1estimates}
Given $\epsilon \in (0, \gamma_{0})$, there exists $M_{\epsilon} \geq 1$
such that for  $ g \in X_{2} \cap \cont, t > 0,$
\begin{equation}
\label{thalf}
\| e^{t \lop_{2}}g \|_{\cont} \leq M_{\epsilon} t^{- \frac{1}{2}}
e^{- \gep t} \| g \|_{\buc} 
\end{equation}
and
\begin{equation}
\| e^{ t \lop_{2}}g \|_{\cont} \leq M_{\epsilon} e^{- \gep t} \|g \|_{\cont},
\label{sgdecay}
\end{equation}
where $\gep = \gamma_{0} - \epsilon$.
\end{lemma}

\proof Lemma \ref{lempartsectorial} implies that the part of $\lop$ in
$\cont$ generates an analytic semigroup in the Banach space $\cont$.
 So there exist $M>0$ and
$\omega \in \rs$ such that for each $t>0, g \in \cont$,
\begin{equation}
\label{semigroup}
\| e^{t \lop}g \|_{\cont} \leq M e^{\omega t} \|g\|_{\cont}.
\end{equation}
 
Fix $\epsilon \in (0, \gamma_{0})$.
We appeal to \cite{Lunardi},  in the notation of which,
 let $\alpha = \half$
and $n=0$. The spaces $D_{\lop} (\half , p), 1 \leq p \leq \infty$ are
defined in \cite[p 45]{Lunardi}; note the last remark on that page. Now
observe that
\begin{equation}
\label{halfone}
D_{\lop}(\half, 1) \hookrightarrow  \buc^{1}.
\end{equation}
This follows from Landau's inequality, \cite[p 46, Proposition 2.2.2  and
p 24, Theorem 1.2.13  with $\theta = \half$]{Lunardi}. This and
\cite[p 59, Proposition 2.3.3  with $\beta= \half$ and $p=1$]{Lunardi}
 together yield the
existence of $\hat{M}>0$ such that
for each
$g \in X_{2} \cap \cont$,
\begin{equation}
\| e^{t \lop_{2}}g\|_{\cont} \leq  \hat{M} t^{- \frac{1}{2}} e^{- \gep t}
\|g\|_{\buc} ~~\mbox{ for each}~ t >0. 
\label{est2}
\end{equation}
In addition,
\begin{equation}
\label{halfbeta}
\buc^{1} \hookrightarrow D_{\lop} (\half, \infty) ~~~\mbox{and}~~~
D_{\lop} (\beta , \infty) \hookrightarrow \buc^{1}, \beta \in (\half, 1),
\end{equation}
by \cite[p 86, Theorem 3.1.12  with $\theta = \half$ and $\theta = \beta$
respectively]{Lunardi}. \cite[p 59, 
Proposition 2.3.3 with $\beta \in (\half, 1), p= \infty$]{Lunardi} and
(\ref{halfbeta}) give the existence of $\hat{M}>0$ such that for
each $g \in X_{2} \cap \cont$,
\begin{eqnarray}
\| e^{t \lop_{2}}g\|_{\cont}  
& \leq & \hat{M}t^{\frac{1}{2} - \beta}e^{-\gep t}\|g\|_{\cont}, ~~\mbox{for
each }~t>0,  \nonumber \\
& \leq & \hat{M} e^{- \gep t} \| g \|_{\cont} ~~\mbox{when}~~ t \geq 1.
\label{beta}
\end{eqnarray}
\noindent It follows from (\ref{semigroup}) and (\ref{beta}) that there exists
$\tilde{M} >0$ such that 
\begin{equation}
\| e^{t \lop_{2}} g \|_{\cont} \leq \tilde{M}
e^{-\gep t}\|g\|_{\cont} ~~\mbox{ for all}~~ t>0.
\label{est}
\end{equation}
(\ref{thalf}) and (\ref{sgdecay}) follow from (\ref{est2}) and (\ref{est}). \qed 

\section{Local stability}

It is useful to  formulate  (\ref{frameeq}) as an abstract ordinary 
differential equation.  Let $T>0$ and let $v \in C (\rs \times
[0,T], \RN)$ be such that $v, v_{t}, v_{x}$ and $v_{xx}$ 
are bounded and uniformly continuous 
on $\rs \times (0,T)$. Define $y(t)(x) = v(x,t)
 - w(x), ~~(x,t) \in \rs \times [0,T],$ where $w$ is the travelling wave
 introduced in
\tw. Then $v$ satisfies (\ref{frameeq}) if  and only if $y \in
C^{1} ((0, T), \buc) \cap C( (0, T), \buc^{2})$ satisfies
\begin{equation}
\label{lopform}
y'(t) = \lop (y(t)) + \rop (y(t)), ~~ t \in (0,T)
\end{equation}
where $\rop : \buc^{1} \rightarrow \buc$ is given by
$$ \rop (y) = f(w+y, w' + y') - f(w , w') -d_{p}f[w,w']y' - d_{q}f[w,w']y, ~~
y \in \buc^{1}.$$
Note that $\rop$ is continuously differentiable, and that
 $\| \rop (y) \|_{\buc} / \|y\|_{\buc^{1}} \rightarrow 0$ as
$\| y \|_{\buc^{1}} \rightarrow 0$. 

Following \cite[p 108]{Henry}, we adopt an elementary approach to
proving local stability, based on the variation of constants formula
and the estimates of Lemma \ref{lemc1estimates}. An alternative is to
use centre-manifold theory and the existence of foliations - see
\cite{BJ}, \cite{BDL}, \cite{CLL}.

\begin{theorem}
\label{localstabthm}
Let $\epsilon \in (0, \gamma_{0})$. Then there exist $\nuep >0, \kep >0$
and $\dep >0$ such that if $\phi \in \buc^{1}$ satisfies
\begin{equation}
\label{c1nbd}
\| \phi - w( \cdot + \chi_{0}) \|_{\buc^{1}} < \nuep
\end{equation}
for some $\chi_{0} \in \rs$, then there exists $\chi_{\infty} \in
[\chi_{0} - \dep , \chi_{0} + \dep]$ such that
\begin{equation}
\label{localconv}
\| v^{\phi} \fixt - w( \cdot + \chi_{\infty}) \|_{\buc^{1}} \leq
\kep e^{- \gep t}, ~~ t>0.
\end{equation}
Note that $ \kep $ and $\dep>0$ are independent of the exact choice
of $\phi , \chi_{0}$ satisfying (\ref{c1nbd}).
\end{theorem}

\proof We  first prove a convergence result for (\ref{lopform}), and then
deduce Theorem \ref{localstabthm} by interpreting this in terms of 
(\ref{frameeq}) and the travelling wave $w$.
  The idea for the proof comes from \cite[p 108, Exercise 6]{Henry}.
For $\chi \in \rs$, define $\what : \rs \rightarrow \buc^{1}$ by
$\what (\chi)(x) = w(x + \chi) - w(x),$ $ x \in \rs$.
Then $\what (0)=0$, and for each $\chi \in \rs$, $\lop \what (\chi ) +
 \rop (\what (\chi)) =0$, since $w(\cdot + \chi)$ is a stationary solution
of (\ref{frameeq}). Since $w$ satisfies \tw~ and $f \in \fcont$,
$w \in \buc^{3}$. So $\what : \rs \rightarrow \buc^{1}$ is twice
 continuously differentiable, and
\begin{equation}
\label{derivative}
d \what [\chi_{0}] \chi = \chi w'( \cdot + \chi_{0}) ~~\mbox{for each}~~
\chi_{0}, \chi \in \rs. 
\end{equation}

Let $\hop (y, \chi) = \wstar ( y - \what (\chi)) \in \rs,
 ~~(y, \chi) \in \buc^{1}
\times \rs$, where $\wstar$ is as in Lemma \ref{lemprojection}. Then
$\hop$ is continuously differentiable, $\hop (0,0) =0$ and 
$d_{\chi}\hop [0,0] \chi =- \chi$ for each $\chi \in \rs$. So it
 follows from the implicit function theorem that there is an
open ball $\ball$ in $\buc^{1}$ 
(centre 0, radius $ \rho_{0}$), an open neighbourhood
$( -\delta_{0} , \delta_{0})$ of $0$ in $\rs$ and a continuously differentiable
function $\zeta : \ball \rightarrow \nbd$ such that $\zeta (0) =0$,
$\hop ( y, \zeta (y))=0$ for $y \in \ball$, and if
$ \hop (y, \chi) =0$ for some 
$y \in \ball, \chi \in \nbd$, then
$ \chi = \zeta (y)$.
By (\ref{projnormal}), we can choose $\rho_{0} >0$ smaller if necessary so 
that
$\wstar (w'( \cdot + \chi )) > \half$ whenever $\chi = \zeta (y)$ 
for some $y \in \ball$.

Proposition \ref{localthm} (Appendix) 
ensures that given initial data $\yz \in 
\buc^{1}$, there is a unique local classical solution $y : 
(0, \tau (\yz)) \rightarrow \buc^{2}$ of (\ref{lopform})
such that $\| y(t) - \yz \|_{\buc^{1}} \rightarrow 0$ as $t \rightarrow 0$.
  For  $\yz \in \ball$, let
$\hat{t} \in ( 0, \tau (\yz))$ be such that $y(t) \in \ball$
for each $t \in [0, \that]$. 
 For such $t$, 
define $\chi (t) = \zeta (y (t))$, where $\zeta$ is as given by
the implicit function theorem above. Then $\chi (t) \in \nbd$ and
$\wstar (y(t)) = \wstar ( \what (\chi (t)))$. Define
$\up (t) = y(t) - \what (\chi (t))$. Since $\wstar (\up (t)) =0$,
$\up (t) \in X_{2}$ (where $X_{2}$ is as defined before Lemma
 \ref{lemprojection}).
Note that $\what (\chi (\cdot)) = \what (\zeta (\yphi(\cdot)))$
 and $\up (\cdot)$
are both continuously differentiable on $(0,\that)$, and 
 since $\yphi \in C^{1}((0, \that), \buc)$ and $X_{2}$ is a 
closed subspace of $\buc$,   $\up'(t) \in X_{2}$
for $0<t<\that$.

Acting on (\ref{lopform}) with $\wstar$ and using (\ref{derivative}),
the fact that $\what (\chi)$ is a stationary solution of (\ref{lopform})
for each $\chi$ and
the properties of $\wstar$ together yield that 
 for $0<t<\that$, 
\begin{equation}
\label{x1proj}
\chi'(t) \wstar (w'(\cdot + \chi (t))) = \wstar (\rop (\what (\chi(t))
+ \up (t)) - \rop (\what (\chi (t)))).
\end{equation}
So
\begin{equation}
\chi '(t) = \Phi (\chi (t), \up (t)), ~t \in (0, \that),
\label{chieqn}
\end{equation} 
where we define
\begin{equation}
\label{defPhi}
\Phi (\chi , \up) = \frac {\wstar (\rop (\what (\chi) + \up) - \rop (\what 
(\chi)))} { \wstar (w' ( \cdot + \chi))}, ~(\chi, \up )
 \in \rs \times \buc^{1}.
\end{equation}
Similarly,  acting on (\ref{lopform})
with $I - \pop$ (see (\ref{contour})) gives that
\begin{equation}
\label{upeq}
\up '(t)  = \lop_{2} \up (t) + \Psi ( \chi (t), \up (t)),~~ t \in (0,
\that),
\end{equation}
where
\begin{equation}
\label{defPsi}
\Psi (\chi , \up ) = (I - \pop) \{ \rop (\what (\chi) + \up )
- \rop (\what (\chi )) \} - (I - \pop) d \what [\chi] 
\Phi (\chi  , \up ) .
\end{equation}

\noindent 
Now for $\up \in \buc^{1}$ and $\chi \in \rs$ with $| \chi | \leq 1$
and  small enough
that $\wstar (w' ( \cdot + \chi)) > \half$,
\begin{eqnarray*}
| \Phi (\chi , \up )|  \leq  2 \| \wstar \|_{\buc^{*}} K(\chi, \up) 
\| \up \|_{\buc^{1}},~~\mbox{where}~~
K(\chi , \up)= \sup_{0 \leq \theta \leq 1}
\{ \| d \rop [ \what (\chi) + \theta \up ]\|_{L(\buc^{1}, \buc)} \}.
\end{eqnarray*}
Since $d \rop [0] =0$,
$K( \chi , \up) \rightarrow 0$ as $ | \chi| +
 \| \up \|_{\buc^{1}} \rightarrow 0$.
Also,
\begin{equation}
\| \Psi (\chi , \up ) \|_{\buc} \leq \| I - \pop \|_{L(\buc , \buc)} 
K( \chi , \up ) \| \up \|_{\buc^{1}} + \| I - \pop \|_{L(\buc , \buc)}
\| d \what [\chi ] \|_{L( \rs , \buc^{1})} |\Phi (\chi , \up )|.
\label{Psiest}
\end{equation}
So, since $\| d \what [ \chi ] \|_{L(\rs , \buc^{1})}$ is bounded independently
of $| \chi | \leq 1$, there exists a constant $\hat{K}>0$ such that
\begin{equation}
|\Phi (\chi , \up)| + \| \Psi ( \chi , \up ) \|_{\buc} \leq \hat{K}
K(\chi , \up) \| \up \|_{\buc^{1}}, ~~\mbox{where}~ K(\chi, \up)
\rightarrow
0 ~\mbox{as}~ |\chi|+ \|\up \|_{\buc^1} \rightarrow 0.
\label{PsiPhiest}
\end{equation}

\noindent Henceforth fix $\epsilon \in (0, \gamma_{0})$. Choose $\sep >0$ so that
\begin{equation}
\mept \sep \int_{0}^{\infty} s^{- \half} e^{-( \gept - \gep) s} \dt
= \mept \sep \int_{0}^{\infty} s^{- \half}e^{- \frac{\epsilon}{2} s} \dt
< \frac{1}{2}, 
\label{sepdef}
\end{equation}
where $\mept \geq 1$ is as in Lemma \ref{lemc1estimates}.
Let $\tilde{K} > 0$ be such that
$K( \chi, \up) < \tk $ whenever $ | \chi | < \delta_{0}$ and
$ \| \up \|_{\buc^{1}} < \rho_{0}$.
Now using (\ref{PsiPhiest}),
 we can choose $\rep \in (0, \rho_{0})$, 
$\dep \in (0, \delta_{0})$ such that $\rep < \frac{\gep^{2}}{2 \hat{K}
 \tilde{K}}$ and
\begin{equation}
\label{repdepconditions}
 \| \Psi (\chi, \up) \|_{\buc}  \leq  \sep \| \up \|_{\buc^{1}}, 
 \| \what (\chi) + \up \|_{\buc^{1}}  \leq  {\displaystyle \frac{\rho_{0}}{2}}
 ~~\mbox{for all}~~ (\chi, \up) ~\mbox{with}~ |\chi|
\leq \dep ~\mbox{and}~ \| \up \|_{\buc^{1}} \leq \rep.
\end{equation}
Let $\nep \in (0, \rho_{0})$ be such that
\begin{equation}
\label{datacondition}
\| \yz \|_{\buc^{1}} < \nep \Rightarrow 
| \zeta(\yz)| <  \dep/2 ~\mbox{and}~ 
\|  \yz \|_{\buc^{1}}+ \| \what (\zeta (\yz)) \|_{\buc^{1}}
 < 
\rep/({2 \mept}).
\end{equation}
Fix initial data $\yz \in \buc^{1}$ with $\| \yz \|_{\buc^{1}} < \nep$. 
Define 
$t_{0} = \sup_{0 \leq t < \tyz} \{ t : y(s) \in \ball ~~\mbox{for all}~~
s \in [0,t] \}$. 
For $t \in [0, t_{0})$, $\chi (t) = \zeta (y(t))$ and $\up (t) = y(t)
- \what (\chi (t))$ are well-defined and have the properties described 
above. 
By the choice of $\nep$, 
$| \chi (0) | < \frac{\dep}{2}$ and $ \| \up (0) \|_{\buc^{1}}
< \frac{\rep}{2 \mept}$.
Define
$m(t) = \sup_{0 \leq s \leq t} \left\{ e^{\gep s} \| \up (s) \|_{\buc^{1}} 
\right\},$ $t \in [0, t_{0})$.
Then since $\up$ satisfies (\ref{upeq}) and $\gep = \gamma_{0} - \epsilon$,
 it follows from the
 variation of constants formula, Lemma \ref{lemc1estimates},
(\ref{sepdef}) and (\ref{repdepconditions}) that 
for $0 \leq s \leq t < t_{0}$, 
\begin{eqnarray*}
e^{\gep s} \| \up (s) \|_{\buc^{1}} & = & 
e^{\gep s} \left\| e^{s \lop_{2}}  \up (0) + \int_{0}^{s} e^{(s- \ts)\lop_{2}}
\Psi (\chi (\ts), \up (\ts)) ~ d \ts \right\|_{\buc^{1}} \\
& \leq & \mept e^{(\gep - \gept)s} \| \up (0)\|_{\buc^{1}}
+ e^{ \gep s} \sep \mept \int_{0}^{s} (s - \ts)^{- \half} e^{- \gept (s- \ts)}
\| \up (\ts) \|_{\buc^{1}} ~ d \ts \\
& \leq & \mept \| \up (0) \|_{\buc^{1}} + \half m(t). 
\end{eqnarray*}
Whence $m(t) \leq 2 \mept \| \up (0) \|_{\buc^{1}}$ for each
 $t \in [0, t_{0})$. It follows, using 
(\ref{chieqn}), (\ref{PsiPhiest}), that
\begin{equation}
\| \up (t) \|_{\buc^{1}} \leq \rep e^{ - \gep t} ~~\mbox{and}~~
| \chi '(t) | = | \Phi (\chi (t) , \up (t)) | \leq \hat{K} \tilde{K}
\rep  e^{ - \gep t}, ~~ t \in (0, t_{0}).
\label{yhatdecay}
\end{equation}
This, together with the facts that $| \chi (0) | < \dep /2$ and  
$\rep < \frac{\gep^{2}}{2 \hat{K}
 \tilde{K}}$,
yields that for each $t \in [0, t_{0})$, 
\begin{equation}
\label{chibound}
| \chi (t) | 
\leq {\dep}/2 + \hat{K} \tilde{K} \rep \gep^{-1}
[1 - e^{- \gep t}] ~ < ~ \dep.
\end{equation}
Now it follows from the definition of $t_0$, (\ref{repdepconditions}),
(\ref{yhatdecay}) and (\ref{chibound}) that $t_0 = \tyz$. 
And  Proposition \ref{condthm} (Appendix)
 shows that if
$\tyz < \infty$, then $\sup_{0 \leq s \leq t} \| y(s) \|_{\buc} \rightarrow
\infty$ as $t \uparrow \tyz$.
So $t_0 = \tyz = \infty$, and  
(\ref{yhatdecay}) and (\ref{chibound}) hold for all
$t \geq 0$. Since $| \chi ' (\cdot ) | \in L^{1} ( (0, \infty), \rs)$ and
$| \chi (t) | \leq \dep $ for all $t \geq 0$, there exists
$\hat{\chi} \in [ -\dep, \dep]$ such that  
\begin{equation}
\label{chiconv}
| \hat{\chi} - \chi (t) | 
\leq \hat{K} \tilde{K} \rep {\gep}^{-1} e^{- \gep t}, ~~ t>0.
\end{equation}
 
We now rewrite (\ref{yhatdecay}) and (\ref{chiconv}) in terms of the
travelling wave $w$. Recall that $y$ is a solution of (\ref{lopform}) with
initial data $\yz$  if and 
only if $v^{\phi}(\cdot ,t) = y(t) + w$ is a solution of (\ref{frameeq})
 with initial
data $\phi = \yz + w$, and that $\what (\chi) (x) = w(x + \chi) - w(x)$
for $x, \chi \in \rs$. 
So $\| \yz \|_{\buc^{1}} =  \| \phi - w\|_{\buc^{1}}$, and  
\begin{equation}
\| \up (t) \|_{\buc^{1}} =   
\| y(t) - \what (\chi (t)) \|_{\buc^{1}} =
\| v^{\phi}( \cdot , t) - w (\cdot + \chi (t)) \|_{\buc^{1}}.
\end{equation}
Hence if $\| \phi - w \|_{\buc^{1}} \leq \nuep$, (\ref{yhatdecay}) and
(\ref{chiconv}) give that
\begin{eqnarray}
\| v^{\phi}(\cdot, t) - w( \cdot + \hat{\chi}) \|_{\buc^{1}} & \leq &
\| v^{\phi}(\cdot, t) - w(\cdot + \chi (t)) \|_{\buc^{1}}
+ \| w( \cdot + \chi (t)) - w (\cdot + \hat{\chi}) \|_{\buc^{1}} \nonumber \\
& \leq & \rep e^{ - \gep t} + | \chi (t) - \hat{\chi} |
\| w'\|_{\buc} \nonumber 
~ \leq ~  K_{\epsilon} e^{- \gep t}, 
\label{finalconv}
\end{eqnarray}
where $K_{\epsilon} = \rep \{ 1 +  \hat{K} \tilde{K} \gep^{-1}
\|w'\|_{\buc}\}
$. 
To complete the proof, note that if $\phi \in \buc^{1}$ satisfies 
$\| \phi - w( \cdot + \chi_{0}) \|_{\buc^{1}} < \nuep$ for some $\chi_{0}
\in \rs$, then $\| \phi (\cdot - \chi_{0} ) - w (\cdot) \|_{\buc^{1}} < \nuep$.
The above analysis immediately implies  that 
$ \| v^{\phi} (\cdot, t) - w(\cdot + \hat{\chi} + \chi_{0} ) \|_{\buc^{1}}
\leq K_{\epsilon} e^{- \gep t} $
for all $t \geq 0$. The result follows. \qed

\section{Global stability for monotone initial data}
We turn now to the global stability of the wave $w$.
Note first that if $\phi \in \buc^{1}$ satisfies $\m \leq \phi(x) \leq \p$
for all $x \in \rs$, then it follows from Theorem \ref{globalthm} (Appendix)
 that the 
initial value problem (\ref{frameeq}) has a unique classical solution
$v^{\phi}$ that exists for all time, and that $\m \leq v^{\phi}(x,t)
 \leq \p$ for
all $x \in \rs, ~t \geq 0$.

  In this section, we consider the initial-value
problem (\ref{frameeq}) with initial data 
 $\phi \in \buc^{1}$  
 satisfying the following
conditions : 
\begin{itemize}
\item[($\phi 1$)] $\phi (x) \rightarrow \epm$ as $x \rightarrow \pm \infty$,
and $\phi '(x) \rightarrow 0$ as $|x| \rightarrow \infty$,
\item[($\phi 2$)] $\phi ' (x) \geq 0$ for each $x \in \rs$. 
\end{itemize}
Our approach is similar to that of \cite[pp 245-248, Theorem 6.1]{Vol3}.

\begin{theorem}
\label{globalmonotonethm}
Let $f \in \fcont$ satisfy {\bf (f1) - (f4)} and $\phi \in \buc^{1}$ 
satisfy {\bf ($\phi 1) - (\phi 2)$}. Then there exists $\chi_{\infty}
\in \rs$ such that for each $\epsilon \in (0, \gamma_{0})$, there exists
$N_{\epsilon}>0$ such that the solution $v^{\phi}$ of (\ref{frameeq})
with initial data $\phi$ satisfies
\begin{equation}
\label{globalmonotonedecay}
\| v^{\phi} (\cdot, t) - w( \cdot + \chi_{\infty}) \|_{\buc^{1}} \leq
N_{\epsilon} e^{- \gep t}, ~~~\mbox{for all}~~t>0.
\end{equation}
\end{theorem}

\proof The idea is to construct a function $\pstar$, from $\phi$ and the wave
$w$, such that the solution $v^{\pstar}$ of (\ref{frameeq}) satisfies
(\ref{globalmonotonedecay}), and then to use a homotopy argument to deduce the
corresponding result for $\phi$. 

Fix $\epsilon \in (0, \gamma_{0})$. We begin with the construction of
$\pstar$. Let $\nuep$ be as in (\ref{datacondition}). Choose $\eta_{1}>0$
sufficiently large that
\begin{equation}
\label{etaz}
\pm x \geq +\eta_{1}  \Rightarrow \| \phi (x) - \epm\|, ~
\|w(x) - \epm \|, ~ \|w'(x) \|, \| \phi '(x) \| < \frac{\nep}{4}.
\end{equation}
Choose $\eta_{2} > \eta_{1}+1$ so that $\phi (\eta_{2}) > w(\eta_{1})$
and $\phi (-\eta_{2}) < w(-\eta_{1})$. Define $\pstar : \rs \rightarrow
\RN$ by $\pstar (x) = w(x)$ for $|x| \leq \eta_{1}$ and $\pstar(x) = \phi (x)$
for $|x| \geq \eta_{2}$; for $|x| \in [\eta_{1}, \eta_{2}]$, define $\pstar
(x)$ so that $\pstar \in \buc^{1}$ is  increasing and $\| (\pstar) '(x)
\| < \nep /4$ for each $x, |x| \geq \eta_{1}$. By construction, 
\begin{equation}
\| \pstar - w \|_{\buc^{1}} < \frac{\nep}{2}.
\label{pstarest}
\end{equation}
Here is the construction that underlies the homotopy argument. As in
\cite[p 246]{Vol3}, define
\begin{equation}
\label{phitdef}
\phi_{\tau}(x) = \min \{ \phi (x), \pstar (x - \tau) \}, ~~ \tau \in \rs,
x \in \rs.
\end{equation}
The minimum is calculated componentwise. For each $\tau$, $\phi_{\tau}$
is clearly continuous and  increasing. It also follows directly from 
(\ref{phitdef}) that for each fixed $x \in \rs$, $\pt (x)$ is a 
decreasing function of $\tau$. 
The following crucial property of $\pt$ is proved in \cite[p 246]{Vol3};
\begin{equation}
\pmin (x) = \phi (x) ~~\mbox{and}~~ \ppos (x) = \pstar (x - 2 \eta_{2})
~\mbox{for all}~ x \in \rs.
\label{pminpos}
\end{equation}

The existence theory for the initial-value problem for (\ref{frameeq}) 
 in the Appendix requires the initial data in  
$\buc^{1}$. We
introduce mollifications of $\pt$ in order to consider $\tau$-dependent 
initial-value problems. For $b \in (0,1)$,
 let $\kb : \rs \rightarrow [0, \infty)$
 be a standard normalised
mollifier, supported in $[-b, b]$ (see, for example, \cite[p 46]{Davies}).
For $\tau \in \rs$, $b \in (0,1), x \in \rs$, let 
\begin{equation}
\label{psitbdef}
\ptb (x) = ( \pt * \kb) (x) = \int_{-\infty}^{\infty} \pt (x-s) \kb (s)~ds.
\end{equation}
By construction, $\m \leq \ptb (x) \leq \p$ for all $x$. It follows from
Theorem \ref{globalthm} (Appendix) 
that the initial-value problem (\ref{frameeq}) with 
initial data $\ptb$ has a unique classical solution $\vptb$ that exists
for all time, and that
\begin{equation}
\label{vptbbound}
\m \leq \vptb (x,t) \leq \p ~~~\mbox{for all}~~ x \in \rs, ~~ t \geq 0.
\end{equation}

The approach is to advance the parameter $\tau$ with step $-h <0$
(to be determined) from $\tau = 2 \eta_{2}$ to $\tau = -2 \eta_{2}$, at
each stage proving that the solution $\vptb$ with initial data $\ptb$
converges in $\buc^{1}$ to a translate of $w$. At $\tau = -2 
\eta_{2}$, the initial data is $\phi * \kb$, by (\ref{pminpos}); letting
$b \rightarrow 0$ will then yield the required result.

We seek $\hep >0$, independent of $b \in (0,1), \tau \in \rs, 
 T \geq 1$, such that
\begin{equation}
\| v^{\psi_{\tau - \hep, b}} (\cdot, T) - v^{\psi_{\tau, b}}(\cdot, 
T )\|_{\buc^{1}} \leq \frac{\nep}{4}.
\label{thatest}
\end{equation}
By Landau's inequality,
\begin{equation}
\| (v^{\psi_{\tau - h, b}}  - 
v^{\psi_{\tau, b}})_{x} (\cdot, T) \|_{\buc} \leq
2 \| ( v^{\psi_{\tau - h, b}}   -  v^{\psi_{\tau, b}})
(\cdot, T)
\|_{\buc}^{\half} 
  \| (v^{\psi_{\tau - h, b}}
- v^{\psi_{\tau, b}})_{xx} (\cdot, T) \|_{\buc}^{\half}
\label{landau}
\end{equation}
for each $b \in (0,1), \tau \in \rs, T \geq 1$ and $h>0$. 
We now  show that the first factor
on the right of (\ref{landau}) is  small when $h$ is small. Note
first that for
$x \in \rs, \tau \in \rs, h>0$,
\begin{equation}
\label{death}
\pt (x) \leq \phi_{\tau -h} (x) \leq \pt (x+ h).
\end{equation}
Since mollification preserves ordering and commutes with translation,
it follows that for  $ b \in (0,1)$,
\begin{equation}
\label{psiorder}
\ptb (x) \leq \psi_{\tau -h, b}(x) \leq \ptb (x+h).
\end{equation}
Now since $f$ satisfies {\bf (f1) - (f2)}, the comparison principle 
Theorem \ref{compthm} (Appendix) yields
that 
\begin{equation}
\label{vorder}
\vptb (x,t) \leq v^{\psi_{\tau -h, b}}(x,t) \leq \vptb (x+h, t), ~~
x \in \rs, t>0.
\end{equation}
So by the Mean Value Inequality, for $t>0, x \in \rs$,
\begin{equation}
\label{mean}
\| v^{\psi_{\tau -h, b}} (x,t) - \vptb (x,t) \| \leq 
\| \vptb (x+h, t) - \vptb (x,t) \| \leq h \| ( \vptb )_{x} \fixt \|_{\buc}.
\end{equation}
By Theorem \ref{boundthm} (Appendix) there exists $K_{1} >0$, independent
of $t \geq 1, \tau \in \rs, b \in (0,1)$, such  $\|( \vptb)_{x} \fixt \|_{\buc}
\leq K_{1}$. Hence for each $h>0$, $ t \geq 1, \tau \in \rs, b \in (0,1)$,
\begin{equation}
\| v^{\psi_{\tau -h, b}} ( \cdot, t) - \vptb (\cdot , t )\|_{\buc}
\leq K_{1} h.
\label{meanvalue}
\end{equation}
It follows from (\ref{vptbbound}) and Theorem \ref{boundthm} that the
second factor on the right of (\ref{landau}) is bounded independently
of $\tau \in \rs, h >0, b \in (0,1), T \geq 1$.
The existence of $\hep >0$ satisfying (\ref{thatest}),
 independent of $b \in (0,1), \tau \in \rs$
and $T \geq 1$,  thus follows from (\ref{landau}) and 
(\ref{meanvalue}). We choose $\hep >0$ smaller if necessary so that there
exists $n \in \N$ such that 
\begin{equation}
\label{integer}
4 \eta_{2} = n \hep.
\end{equation}

Now $\| \pstar * \kb - w \|_{\buc^{1}} \rightarrow 0$ as $b \rightarrow 0$,
so it follows from (\ref{pstarest}) that for $b \in (0, b_{0})$ say,
$\| \pstar * \kb - w \|_{\buc^{1}} < \nuep$. Hence
$\| \psi_{2 \eta_{2}, b} - w(\cdot - 2 \eta_{2}) \|_{\buc^{1}} < \nuep$.
With $\gep, K_{\epsilon}, \dep >0$ (independent of $b$) as in Theorem
\ref{localstabthm},  there
exists $\chi_{2 \eta_{2}, b} \in [ -2 \eta_{2} - \dep, -2 \eta_{2} + \dep]$ 
such that
\begin{equation}
\label{one}
\| v^{\psi_{2 \eta_{2}, b}} \fixt - w(\cdot + \chi_{2 \eta_{2}, b}) \|_{\buc^{1}}
\leq \kep e^{- \gep t} ~~~\mbox{for all}~~ t >0.
\end{equation}
Next define
\begin{equation}
\label{tep}
\tep = \max \{ 1 ,\frac{1}{\gep}\log \frac{4 \kep}{\nep} \}.
\end{equation}
(Clearly $\tep$ is independent of $b \in (0,b_{0})$.) So by (\ref{one})
and (\ref{tep}),
\begin{equation}
\label{two}
\| v^{\psi_{2 \eta_{2}, b}} (\cdot, \tep) - w (\cdot + \chi_{2 \eta_{2}, b})
\|_{\buc^{1}} \leq \frac{\nuep}{4}.
\end{equation}
Together with (\ref{thatest}), this yields that
\begin{equation}
\label{four}
\| v^{\psi_{2 \eta_{2} - \hep, b}} (\cdot, \tep) - w( \cdot + \chi_{2 \eta_{2},
b}) \|_{\buc^{1}} \leq \frac{\nuep}{2}.
\end{equation}
So by Theorem \ref{localstabthm}, there exists $\chi_{2 \eta_{2} - \hep, b}
\in [\chi_{2 \eta_{2}, b} - \dep, \chi_{2 \eta_{2}, b} + \dep ]
\subset [-2 \eta_{2} -2 \dep, -2 \eta_{2} + 2 \dep]$ such that for $t > \tep$,
\begin{equation}
\label{five}
\| v^{\psi_{2 \eta_{2} - \hep, b}} \fixt - w( \cdot + \chi_{2 \eta_{2} - \hep,
b} ) \|_{\buc^{1}} \leq \kep e^{- \gep ( t - \tep)}.
\end{equation}
Arguing by induction, it follows that given $m \in \N$, there exist
$\chi_{2 \eta_{2} - k \hep, b} \in [-2 \eta_{2} -k \dep, -2 \eta_{2} + k \dep]$
for each $0 \leq k \leq m$ such that
\begin{equation}
\label{six}
\| v^{\psi_{2 \eta_{2} - m \hep, b}} (\cdot, m \tep) - w( \cdot + 
\chi_{2 \eta_{2} - (m-1) \hep, b}) \|_{\buc^{1}} \leq \frac{\nep}{2},
\end{equation}
and for $t > m \tep$,
\begin{equation}
\label{seven}
\| v^{\psi_{2 \eta_{2} - m \hep, b}} \fixt - w( \cdot+
\chi_{ 2 \eta_{2} - m \hep, b}) \|_{\buc^{1}} \leq \kep e^{- \gep (t - m \tep)}.
\end{equation}

In particular, (\ref{six}) and (\ref{seven}) hold for $n$ satisfying 
(\ref{integer}). Since $\psi_{-2 \eta_{2}, b} = \phi * \kb$, this yields that
for each $b \in (0, b_{0})$, there exists 
 $\chi_{-2 \eta_{2} + \hep , b} \in [-2 \eta_{2} - (n-1) \dep,
-2 \eta_{2} + (n-1) \dep]$ such that
\begin{equation}
\| v^{\phi * \kb } ( \cdot, n \tep) - w (\cdot + \chi_{-2 \eta_{2} + \hep, b})
\|_{\buc^{1}} \leq \frac{\nep}{2}.
\label{eight}
\end{equation}
We now let $b \rightarrow 0$. 
The interval $[ -2 \eta_{2} - (n-1) \dep, -2 \eta_{2} + (n-1) \dep]$ is
independent of $b \in (0, b_{0})$. So there is a sequence $\{ b_{k} \}
\subset (0, b_{0}), b_{k} \downarrow 0$ and $\chi_{\epsilon} \in
[ -2 \eta_{2} - (n-1) \dep , -2 \eta_{2} + (n-1) \dep]$ such that
\begin{equation}
\label{nine}
\chi_{-2 \eta_{2} + \hep, b_{k}} \rightarrow \chi_{\epsilon} ~~\mbox{as}~~
k \rightarrow \infty.
\end{equation}
Thus there exists $k_{0} \in \N$ such that 
\begin{equation}
\label{ten}
k \geq k_{0} \Rightarrow \| w( \cdot + \chi_{-2 \eta_{2} + \hep, b_{k}})
- w(\cdot + \chi_{\epsilon}) \|_{\buc^{1}} \leq \frac{\nuep}{4}.
\end{equation}
Proposition \ref{localthm} (Appendix) yields the existence of $r, K >0$ such that
for $n$ as in (\ref{integer}) and $\hat{\phi}, \tilde{\phi} \in \buc^{1}$,
\begin{equation}
\label{continitial}
\| \hat{\phi} - \tilde{\phi} \|_{\buc^{1}} \leq r \Rightarrow 
\| v^{\hat{\phi}} (\cdot , n \tep ) - v^{\tilde{\phi}} (\cdot , n \tep)
\|_{\buc^{1}} \leq K \| \hat{\phi} - \tilde{\phi} \|_{\buc^{1}}.
\end{equation}
Hence since $\| \phi - \phi * \kb \|_{\buc^{1}} \rightarrow 0$ as $ s 
\rightarrow 0$, there exists $k_{1} \in \N$ such that 
\begin{equation}
\label{eleven}
k \geq k_{1} \Rightarrow \| v^{\phi}( \cdot, n \tep ) -
v^{\phi * \kappa_{b_{k}}}( \cdot , n \tep) \|_{\buc^{1}} \leq \frac{\nep}{4}.
\end{equation}
So by (\ref{eight}), (\ref{ten}) and (\ref{eleven}),
\begin{equation}
\label{twelve}
\| v^{\phi}( \cdot , n \tep) - w ( \cdot + \chi_{\epsilon}) \|_{\buc^{1}}
\leq \nep.
\end{equation}
Theorem \ref{localstabthm} yields that
\begin{equation}
\label{thirteen}
\| v^{\phi} \fixt - w (\cdot + \chi_{\epsilon}) \|_{\buc^{1}} \leq
\kep e^{- \gep (t - n \tep)} ~~\mbox{for}~~ t > n \tep.
\end{equation}
Since $v^{\phi}$ is independent of $\epsilon$, and $w$ is not periodic, 
it is immediate that $\chi_{\epsilon_{1}} = \chi_{\epsilon_{2}}$ for
any $\epsilon_{1} , \epsilon_{2} \in (0, \gamma_{0})$.
The result follows. \qed

\section{Global stability for general initial data}

We will invoke an idea from \cite{RTV}. First a preliminary lemma, 
which is a modification of \cite[Lemma 3.3]{RTV}. This result will
be used later, in the proof of Theorem \ref{increasingthm}, as part of
an argument by contradiction.
\begin{lemma}
\label{harnacklem}
Let $D, G : \rs \times [0, \infty) \rightarrow M^{N \times N}$
be continuous $N \times N$-matrix-valued-functions, uniformly
bounded on $ \rs \times [0, \infty)$, such that $D \xt$ is diagonal
and the off-diagonal elements of $G \xt$ are non-negative for each
$ \xt \in \rs \times [0, \infty)$. Let $\ups$ be a non-negative,
uniformly bounded
solution of 
\begin{equation}
\label{linearups}
\ups_{t} \xt = A \ups_{xx} \xt + D \xt \ups_{x} \xt + G \xt \ups \xt,
~~ \xt \in \rs \times [0, \infty),
\end{equation}
such that $\ups_{t} $ is uniformly bounded for $t \geq \half$ and there exist
$\mu_{0}, M_{0} >0$ such that for each $t \geq 0$, 
\begin{equation}
\label{supmax}
\sup_{x \in \rs} \left( \min_{1 \leq i \leq N} \ups_{i} \xt \right) = 
\max_{|x| \leq M_{0} } \left( \min_{1 \leq i \leq N} \ups_{i} \xt
\right) \geq \mu_{0}.
\end{equation}
Then for each $M \geq M_{0}$, there exists $\alpha(M) >0$ such that
for all $t \geq 1$,
\begin{equation}
\label{minmin}
\min_{|x| \leq M} \min_{1 \leq i \leq N} \ups_{i} \xt \geq \alpha(M).
\end{equation}
\end{lemma}

\proof Let $M \geq M_{0}$ and recall that $\eg = (1, \ldots, 1)$. It follows
from (\ref{supmax}) that for each $T \geq 0$, there exists $x^{T} \in 
[-M_{0}, M_{0}]$ such that $\ups ( x^{T}, T) \geq \mu_{0} \eg$.
Furthermore, $\ups_{t} \xt$ is bounded independently of $x \in \rs,
 t \geq \half$,
so there exists $T_{0} \in (0, \half)$, independent of $ T \geq 1$, such that
\begin{equation}
T \geq 1, | \that | \leq T_{0} \Rightarrow \ups (x^{T}, T+ \that) \geq
\frac{\mu_{0}}{2} \eg.
\label{utbound}
\end{equation}
We will construct a strictly positive function which lies beneath
$\ups \xt$ for all $t \geq 1$.
By the hypotheses on $D$ and $G$, there are constant
diagonal matrices $D^{-}, D^{+}$ and a constant negative-definite
diagonal matrix $G^{-}$ such that
\begin{equation}
\label{GDdominate}
G^{-}_{ij} \leq G_{ij} \xt ~~\mbox{and}~~~ D^{-}_{ij} \leq D_{ij} \xt \leq
D^{+}_{ij} ~~\mbox{for all}~~ x \in \rs, t \geq 0, i,j \in \{1, \ldots, N\}.
\end{equation}
Consider the two initial-boundary-value problems for
$\upos : [0, \infty) \times [0, 2T_{0}] \rightarrow \RN$ and
$\umin : (-\infty, 0] \times [0, 2T_{0}] \rightarrow \RN$;
$$
\begin{array}{l}
u^{\pm}_{t} = A u^{\pm}_{xx} + D^{\pm} u^{\pm}_{x} + G^{-}u^{\pm},~~  
\pmxt \in (0, \infty) \times (0, 2T_{0}), \\
u^{\pm} (0, t) = \frac{\mu_{0}}{2} \eg ~\mbox{for}~ t \in [0, 2T_{0}],\\
 u^{\pm} (x,0)=0 ~\mbox{for}~ \pm x \in (0, \infty), \\
u^{\pm} \xt \rightarrow 0 ~\mbox{as}~\pm x \rightarrow \infty.
\end{array}
$$
Since $A, D^{\pm}, G^{-}$ are diagonal, we can solve these explicitly
using Laplace transforms to find that for each $i \in \{ 1, \ldots, N\}$
and
$ \pmxt \in (0, \infty) \times
[0, 2T_{0}]$,
\begin{equation}
u^{\pm}_{i} \xt = \pm  \left( \frac{A_{i}}{4 \pi} \right)^{\half} x
e^{-\half D_{i}^{\pm}x} \int_{0}^{t} s^{-\frac{3}{2}} \exp \left[
- \left( \frac{ ( D_{i}^{\pm})^{2}}{A_{i}} - G_{i}^{-} \right) s -
\frac{ A_{i} x^{2}}{4 s} \right] ~ d s. 
\label{uzform}
\end{equation}
We will show that $u^{+}_{x} \xt <0$ for all $x >0, t>0$.
 (\ref{uzform}) yields that for each $i \in \{ 1, \ldots, N\}$,
$$
\left( u^{+}_{i} \right)_{x} \xt =
 \left( \frac{A_{i}}{4 \pi} \right)^{\half} 
e^{-\half D_{i}^{+}x} \int_{0}^{t}
\left\{ 1 - \frac{D^{+}_{i} x}{2} - \frac{A_{i} x^{2}}{2 s} \right\}
s^{-\frac{3}{2}} \exp \left[
- \left( \frac{ ( D_{i}^{+})^{2}}{A_{i}} - G_{i}^{-} \right) s -
\frac{ A_{i} x^{2}}{4 s} \right] ~ d s.
$$
Fix $t \in (0, 2T_{0}]$ and let
 $$x^{t}_{+} = \inf \{ x >0 : u^{+}_{x} (s,t) <0 ~~\mbox{for each}~~
s \in [x, \infty) \}. $$
  The formula for $u^{+}_{x}$ shows that  $u^{+}_{x} \xt
<0$ for $x$ sufficiently large. So $x^{t}_{+} \in [0, \infty)$. Suppose
that $x^{t}_{+}>0$. Then for some $i \in \{ 1, \ldots, N\}$, 
$(u^{+}_{i})_{x} (x^{t}_{+}, t) = 0$. So
$$ A_{i} (u^{+}_{i})_{xx} (x^{t}_{+}, t) = (u^{+}_{i})_{t}
(x^{t}_{+}, t) - G^{-}_{ii}u^{+}_{i} (x^{t}_{+}, t). $$
Now $u^{+}_{i} >0$, $G^{-}_{ii}<0$ and it is clear from
(\ref{uzform}) that $(u^{+}_{i})_{t} >0$. So since
$A_{i}>0$, $(u^{+}_{i})_{xx} (x^{t}_{+},t) >0$.
But this implies that $u^{+}_{i}$ has a strict local minimum at 
$x^{t}_{+}$, which contradicts the fact that $(u^{+}_{i} )_{x} \xt <0$
for all $x > x^{t}_{+}$. Whence $x^{t}_{+}=0$. A similar argument
shows that $u^{-}_{x} \xt >0$ for each $t>0, x <0$.\\[1mm]

\noindent Fix $T \geq 1$. Let $\utpos : [x^{T}, \infty) \times
[T-T_{0}, T+T_{0}] \rightarrow \RN$, $\utmin : (-\infty, x^{T}] \times
[T-T_{0}, T+T_{0}] \rightarrow \RN$ denote the unique solutions of the
two
initial-boundary-value problems
\begin{equation}
\label{utvt}
\begin{array}{l}
u^{T, \pm}_{t} = A u^{T, \pm}_{xx} + D^{\pm} u^{T, \pm}_{x} + G^{-}u^{T, \pm},
~~
( \pm \{x-x^{T}\}, t) \in (0, \infty) \times (T-T_{0}, T+T_{0}), \\
u^{T, \pm} (x^{T}, t) = \frac{\mu_{0}}{2} \eg ~\mbox{for}~
 t \in [T-T_{0}, T+T_{0}], \\
 u^{T, \pm} (x,T-T_{0})=0 ~\mbox{for}~ \pm \{x - x^{T}\} \in (0, \infty), \\
u^{T, \pm} \xt \rightarrow 0 ~\mbox{as}~ \pm x \rightarrow \infty.
\end{array}
\end{equation}
Clearly,
$$
u^{T, \pm} \xt = u^{\pm} (x-x^{T}, t-T+T_{0}), ~~~ (\pm \{x-x^{T}]\}, t)
 \in [0, \infty) \times
[T-T_{0}, T+T_{0}]. $$
So, since $ \pm \upm_{x} <0$ for $t, \pm x >0$, 
\begin{eqnarray}
\min_{x \in [x^{T}, M]} u^{T, +} (x,T) \geq \min_{x \in [0, M_{0}+M]}
u^{+} (x, T_{0})& =& u^{+}(M_{0} + M, T_{0}), \label{uTbound}\\
\min_{x \in [-M, x^{T} ]} u^{T, -} (x,T) \geq \min_{x \in [-M_{0}-M, 0]}
u^{-} (x, T_{0})& =& u^{-} (-M_{0}-M, T_{0}). \label{vTbound}
\end{eqnarray}

Now $u^{T, \pm} \xt >0, \pm  u^{T, \pm}_{x} \xt <0 $ for  $( \pm \{x-x^{T}\},
 t)
 \in
 (0, \infty) \times (T-T_{0}, T+T_{0})$,
so it follows from (\ref{GDdominate}) and (\ref{utvt}) that for such $\xt$,
\begin{equation}
\label{UTsub}
u^{T, \pm}_{t} \xt - A u^{T, \pm}_{xx} \xt
 - D \xt u^{T, \pm}_{x} \xt - G \xt u^{T, \pm} \xt \leq 0.
\end{equation}
So since (\ref{utbound}) holds and $\ups$ is non-negative,
 it follows from the positivity theorem Theorem \ref{positivethm} (i)
(Appendix) that
\begin{equation}
\label{hbound}
\ups \xt \geq u^{T, \pm} \xt, ~~~ (\pm \{x-x^{T}\}, t)
 \in [0, \infty) \times
[T-T_{0}, T+T_{0}]. 
\end{equation}
Hence by (\ref{uTbound}), (\ref{vTbound}), (\ref{hbound}),
\begin{equation}
\label{upslowerbound}
\min_{x \in [-M, M]} \ups (x, T) \geq \min \{ u^{0} (M_{0} +M), T_{0},
v^{0}(-M_{0}-M, T_{0}) \}.
\end{equation}
The right-hand side of (\ref{upslowerbound}) is a strictly positive number
independent of $x \in [-M, M], T \geq 1$. The result follows. \qed

\noindent For $\phi \in \buc^{1}$, define its omega limit set
\begin{equation}
\label{Wdef}
W(\phi) = \{ \psi \in \buc^{1} : ~\mbox{there is a sequence}~
t_{n} \rightarrow \infty ~\mbox{such that}~ \| v^{\phi} (\cdot, t_{n})
- \psi \|_{\buc^{1}} \rightarrow 0 \}.
\end{equation}
Theorem \ref{frontglobalthm} (Appendix) gives conditions on the initial
data $\phi$ under which wave-dependent sub- and super-solutions for
 (\ref{frameeq}) can be constructed. This yields important information
about $W(\phi)$. 

\begin{lemma}
\label{Wproplem}
Let $\hat{\eta} >0$ be as in Theorem \ref{frontglobalthm} (Appendix),
 and let $\phi
\in \buc^{1}$
satisfy (\ref{trap}), (\ref{limsup}) for some $\eta \in (0, \hat{\eta})$.
Then
\begin{itemize}
\item[(i)] $W(\phi)$ is nonempty and compact;
\item[(ii)] there exists $\hat{x}(\phi) \in \rs$
such that for all $x \in \rs, \psi \in W(\phi)$,
$$ w(x-\hat{x}(\phi)) \leq \psi (x) \leq 
w(x +  \hat{x}(\phi)); $$
\item[(iii)] $(\psi)' (x) \rightarrow 0$ as $|x| \rightarrow \infty$
for each $\psi \in W(\phi)$;
\item[(iv)] if $\psi \in W(\phi)$, then $v^{\psi} \fixt \in W(\phi)$
for all $t \geq 0$ and $W(\psi)
 \subset W(\phi)$.
\end{itemize}
\end{lemma}

\proof The {\em a priori} estimates of Theorem \ref{boundthm} (Appendix), the
 Arzela-Ascoli
theorem and estimate (\ref{timetrap}) of Theorem \ref{frontglobalthm}
 (Appendix) together
show  (i). Estimate (\ref{timetrap}) also yields  (ii).  (iii)
follows from (ii), Theorem \ref{boundthm} 
 and Landau's inequality on
a half-line.
(iv) is  a consequence of definition (\ref{Wdef}), the last part of
Proposition \ref{localthm} (Appendix)
 and the semigroup property of solutions of 
(\ref{frameeq}).
\qed

The next theorem is the key. We include a proof for completeness; the
approach is a  minor modification of \cite[Lemma 3.4]{RTV}.

\begin{theorem}
\label{increasingthm}
Let $\phi \in \buc^{1}$ be as in Lemma \ref{Wproplem}. Then there exists
$\psi_{0} \in W(\phi)$,   $\psi_{0}
 (x) \rightarrow \epm$
as $x \rightarrow \pm \infty$, $\left( \psi_{0} \right)
'(x) \rightarrow 0$ as
$|x| \rightarrow \infty$ and $\left( \psi_{0} \right) 
'(x) \geq 0$ for each $x \in \rs$.
\end{theorem}

\proof Define $\fop : W(\phi) \rightarrow [0, \infty]$ by
\begin{equation}
\fop (\psi) = \inf \{ \chi_{0} >0 : \psi (x + \chi)
\geq \psi (x) ~\mbox{for all}~ \chi \geq \chi_{0}, x \in \rs\}.
\end{equation}
Note that since $W(\phi) \subset \buc$, $\psi (x + \fop (\psi)) \geq
\psi (x)$ for each $x \in \rs, ~\psi \in W(\phi)$. 
Lemma \ref{Wproplem}  (ii) shows that $\fop (\psi) < \infty$
for each $\psi \in W(\phi)$. It follows from Lemma \ref{Wproplem}
 (i) that $\fop$ attains its minimum $\fop_{0}$ at a point $\pdo \in 
W(\phi)$. Lemma \ref{Wproplem}  (ii), (iii) ensure that $\pdo (x) 
\rightarrow \epm$ as $x \rightarrow \pm \infty$ and $(\pdo) ' (x) \rightarrow
0$ as $|x| \rightarrow \infty$.

If $\fop_{0} =0$, then $(\pdo) '(x) \geq 0$ for each $x \in \rs$.
So suppose, for contradiction, that $\fop_{0} >0$. We consider the
solution $\vpdo$ of (\ref{frameeq}) with initial data $\pdo$. 
Note first that Lemma \ref{Wproplem} (iv) states that $\vpdo \fixt
\in W(\phi)$ for all $t \geq 0$. By the choice of
 $\pdo$ as the minimiser of $\fop$
 and Theorem \ref{compthm} (Appendix), 
$\fop (\vpdo \fixt) = \fop_{0}$ for all $t \geq 0$, so
$$ \vpdo (x + \fop_{0} , t) \geq \vpdo \xt ~~\mbox{for all}~~x \in \rs,
 t \geq 0.$$ In fact, since $\fop_{0} >0$ and 
Lemma \ref{Wproplem} { (ii)} holds, there exist $\mu_{0}>0,
M_{0} >0$ such that for all $t \geq 0$,
\begin{equation}
\label{muzero}
\sup_{x \in \rs} \left( \min_{1 \leq i \leq N} \vpdo_{i} (x+ \fop_{0}, t)
- \vpdo_{i} \xt \right) = \max_{|x| \leq M_{0}} \left( \min_{1 \leq i \leq N}
\vpdo_{i} (x + \fop_{0}, t) - \vpdo_{i} \xt \right) \geq \mu_{0}.
\end{equation}
Let $\qz, \eepm, \nu$ be as in the preamble to Theorem \ref{fronttrapthm}
(Appendix).
 By Theorem \ref{boundthm} (Appendix), $\| \vpdo_{xx} \fixt \|_{\buc}$
is bounded independently of $t \geq 1$. So it follows from Lemma \ref{Wproplem}
(ii) and Landau's inequality on a half line that $\vpdo_{x} \xt \rightarrow 0$
as $|x| \rightarrow \infty$ at a rate independent of $t \geq 1$. Thus Lemma
\ref{Wproplem} (ii) and (\ref{fprop}) give that there exists
 $M \geq M_{0}$ such that for all $t \geq 0$, 
$\sigma \in [0,1]$ and each
$\fop \in [0, \fop_{0}]$,
\begin{eqnarray}
\pm x &\geq& M  \Rightarrow \nonumber \\ & & d_{q} f[ \sigma \vpdo (x+ \fop,
 t) + (1- \sigma) \vpdo \xt , \sigma \vpdo_{x} (x + \fop,
 t) + (1- \sigma) \vpdo_{x} \xt ]  \eepm
\leq - \frac{\nu}{2} \eepm. 
\label{mplus} 
\end{eqnarray}
For $\delta \geq 0$, define
\begin{equation}
\label{hgoth}
 \dsold \xt = \vpdo (x+ \fop_{0} - \delta, t) - \vpdo \xt, ~~ x \in \rs, 
t \geq 0. 
\end{equation}
Then $\dsol^{0} \geq 0$, and for each $\delta \geq 0$, $\dsold$ is a solution of
\begin{equation}
\label{dequation}
\dsold_{t} \xt = A \dsold_{xx} \xt + c \dsold_{x} \xt + D^{\delta} \xt
\dsold_{x} \xt + G^{\delta} \xt \dsold \xt,
\end{equation}
where
$$
D^{\delta} \xt  =  \int_{0}^{1} d_{p} f[ \sigma \vpdo (x+ \fop_{0}
-\delta, t) + (1- \sigma) \vpdo \xt , \sigma \vpdo_{x} (x + \fop_{0} -
\delta , t) + (1- \sigma) \vpdo_{x} \xt ] ~ d \sigma , $$
$$
G^{\delta} \xt  =  \int_{0}^{1} d_{q} f[ \sigma \vpdo (x+ \fop_{0}
-\delta, t) + (1- \sigma) \vpdo \xt , \sigma \vpdo_{x} (x + \fop_{0} -
\delta , t) + (1- \sigma) \vpdo_{x} \xt ] ~ d \sigma.
$$
Since $f$ satisfies { (f1)} and { (f2)}, the matrices $cI+D^{0}, G^{0}$
satisfy the hypotheses on $D, G$ respectively in   Lemma \ref{harnacklem}.
Also, Theorem \ref{boundthm} (Appendix) shows that $\dsol^{0}_{t} \xt$
is bounded independently of $x \in \rs, t \geq \half$. 
 So with $M$ as in (\ref{mplus}),  Lemma \ref{harnacklem}  (applied to
the function $\dsol^{0}$)  together with
 (\ref{muzero}) imply the existence
of $\alpha (M) >0$ such that for each $t \geq 1$,
\begin{equation}
\label{alphaM}
|x| \leq M \Rightarrow \dsol^{0} \xt \geq \alpha (M).
\end{equation}
Theorem \ref{boundthm}  
shows that $(\vpdo)_{x} \xt$ is bounded independently
of $x \in \rs, t \geq 1$. So there exists $\delta(M) \in (0, \fop_{0})$
such that for each $t \geq 1$,
\begin{equation}
\label{halfalphaM}
|x| \leq M, \delta \in [0, \delta(M) ] \Rightarrow \dsold \xt \geq \half
\alpha(M).
\end{equation}
Now by Lemma \ref{Wproplem} { (i)} and { (iv)}, there is a 
sequence $t_{n} \rightarrow \infty$ and $\pddo \in W(\pdo) \subset W(\phi)$
such that
\begin{equation}
\label{vpdoconv}
\| \vpdo ( \cdot, t_{n}) - \pddo \|_{\buc^{1}} \rightarrow 0 ~~\mbox{as}~~
n \rightarrow \infty.
\end{equation}

Fix $\delta \in [0, \delta(M)]$.
(\ref{halfalphaM}) and (\ref{vpdoconv}) show that $\pddo (x + \fop_{0}
-\delta) \geq \pddo (x)$ for all $x \in [-M, M]$. 
Consider $x \geq M$.  Now 
 $\dsol^{\delta}$ is uniformly bounded, $\dsol^{\delta}
 (M, t) \geq \half \alpha(M)$ for $t \geq 1$ and (\ref{mplus}) holds.
 So Theorem \ref{positivethm} (i)
 (Appendix) (applied to (\ref{dequation})) 
shows that there is a constant $K^{\delta}>0$
 such that for all $x \geq M, t \geq 1$,
\begin{equation}
\dsol^{\delta} \xt \geq - K^{\delta}  e^{- \frac{\nu}{2} t} \eep.
\end{equation}
Whence $\pddo (x + \fop_{0} - \delta) \geq \pddo (x)$ for each $x \geq M$.
Similarly, $\pddo (x+ \fop_{0} - \delta) \geq \pddo (x)$ for $x \leq -M$.
So
\begin{equation}
\label{contradict}
\pddo (x + \fop_{0} - \delta) \geq \pddo (x) ~~\mbox{for all}~~ x \in \rs,
~\delta \in [0, \delta(M)].
\end{equation}
But it follows from (\ref{vpdoconv}) and the fact that $\fop (\vpdo \fixt) 
= \fop_{0}$ for all $t \geq 0$ that $\fop (\pddo) = \fop_{0}$. This
contradicts (\ref{contradict}). Thus $\fop_{0} = 0$ and the result follows.
\qed

\noindent The main result of this paper is the following.

\begin{theorem}
\label{mainthm}
Let $f \in \fcont$ satisfy {\bf (f1) - (f4)}. Let $\hat{\eta} >0$ be as
 in Theorem \ref{frontglobalthm}, and let $\phi$ 
satisfy (\ref{trap}), (\ref{limsup}) for some $\eta \in (0, \hat{\eta})$.
Then there exists $\chi_{\infty}
\in \rs$ such that for each $\epsilon \in (0, \gamma_{0})$, there exists
$N_{\epsilon}>0$ such that the solution $v^{\phi}$ of (\ref{frameeq})
with initial data $\phi$ satisfies
\begin{equation}
\| v^{\phi} (\cdot, t) - w( \cdot + \chi_{\infty}) \|_{\buc^{1}} \leq
N_{\epsilon} e^{- \gep t}, ~~~\mbox{for all}~~t>0.
\end{equation}
\end{theorem}

\proof Theorem \ref{increasingthm}, Theorem \ref{globalmonotonethm}
and Lemma \ref{Wproplem} { (iv)} show that there exists $\chi_{\infty}
\in \rs$ such that $w(\cdot + \chi_{\infty}) \in  W(\phi)$. The result then
follows from Theorem \ref{localstabthm}. \qed

\noindent This implies a uniqueness result for travelling-wave solutions
of (\ref{origeq}).

\begin{corollary}
Let $f \in \fcont$ satisfy {\bf (f1)-(f4)}.  Let $\hat{\eta} >0$ be as
 in Theorem \ref{frontglobalthm}, and let $\phi \in \buc^{1}$ 
satisfy (\ref{trap}), (\ref{limsup}) for some $\eta \in (0, \hat{\eta})$.
Suppose that there exists $\hat{c} \in \rs$ such that $u \xt : = 
\phi (x - \hat{c} t)$ is a travelling-wave solution of (\ref{origeq}).
Then $\hat{c} = c$ and there exists $\chi_{\infty} \in \rs$ such that
$\phi (\cdot) = w(\cdot + \chi_{\infty})$. (Here $w, c$ are as in \tw.)
\label{uniquecor}
\end{corollary}

\proof Theorem \ref{mainthm} shows that there exists $\chi_{\infty} \in \rs$
such that 
\begin{equation}
\label{twowaves}
\| u( \cdot +ct , t) - w( \cdot + \chi_{\infty} ) \|_{\buc^{1}} = 
\| \phi (\cdot + \{ c - \hat{c} \} t) - w( \cdot + \chi_{\infty}) \|_{\buc^{1}}
\rightarrow 0 ~~\mbox{as}~~ t \rightarrow \infty.
\end{equation}
Suppose that $c > \hat{c}$. Since $w(x) \rightarrow \m$ as $x \rightarrow 
-\infty$, we can choose $\hat{x} \in \rs$ such that $w(\hat{x} + \chi_{\infty})
< \p - \hat{\eta} \eep$. But since $\phi$ satisfies (\ref{limsup})
and $c - \hat{c} >0$, $\phi(\hat{x} + \{ c-\hat{c} \} t) > \p - \hat{\eta} 
\eep$
for $t$ sufficiently large. This contradicts (\ref{twowaves}), so
 $c \leq \hat{c}$. A similar argument shows that $c \geq \hat{c}$. Whence
$c= \hat{c}$. The result now follows from (\ref{twowaves}). \qed

\appendix
\section{Appendix}
\vspace{3mm}
{\bf  Comparison theorem}\\[2mm]
For $T>0$, define
$$ \Gamma_{T} = \{v \in C(\rs  \times [0,T], \RN) :
v_t, v_x, v_{xx} ~\mbox{are continuous on}~ \rs 
\times (0,T) \}, $$ 
$$  \Gamma_{T}^{+} = \{v \in C([0, \infty)  \times [0,T], \RN) :
v_t, v_x, v_{xx} ~\mbox{are continuous on}~ (0, \infty) 
\times (0,T) \}. $$
For $v \in \Gamma_{T}$, $(x,t) \in \rs \times (0, T]$ (or
 $v \in \Gamma_{T}^{+}$, $(x,t) \in (0, \infty) \times (0, T)$),  define
\begin{equation}
\mop (v)(x,t) = - v_{t} (x,t) + A v_{xx}(x,t) + D(x,t) v_{x}(x,t) + G(x,t)
v(x,t),
\label{defM}
\end{equation}
and
\begin{equation}
\nop (v)(x,t) = -v_t(x,t) + Av_{xx}(x,t) +c v_{x}(x,t)+ f(v(x,t), v_x(x,t)),
\label{defN}
\end{equation}
where $A$ satisfies {\bf (a)}, $c \in \rs$, 
$f \in \fcont$ satisfies {\bf (f1) - (f2)}
and $D, G : \rs \times [0,T] \rightarrow M^{N \times N}$ are continuous
$N \times N$ matrix-valued functions, bounded on $\rs \times [0,T]$, such that
$D$ is diagonal and the off-diagonal elements of $G$ are non-negative. 
\cite[p 241, Lemma 5.2 and p 242, Theorem 5.3]{Vol3}
yield the following positivity results.
\begin{theorem}
\label{positivethm}
\begin{itemize}
\item[(i)] Let $v \in \Gamma_{T}^{+}$ be such that $v$ is bounded
on $[0, \infty) \times [0, T]$ and $\mop (v) \xt \leq 0$
for $ \xt \in (0, \infty) \times (0, T]$. If $v(x,0) \geq 0$ for all 
$x \in \rs$ and $v(0, t) \geq 0$ for each $t \in [0, T]$, then
$v \xt \geq 0$ for all $\xt \in [0, \infty) \times [0, T]$.
\item[(ii)] 
Let $v \in \Gamma_{T}$ be such that $v$ is bounded on $\rs \times [0, T]$
and $\mop (v)(x,t) \leq 0$ for 
$(x,t) \in \rs \times (0,T]$. If $v(x,0) \geq 0$ for all $x \in \rs$,
then $v (x,t) \geq 0$ for all $(x,t) \in \rs \times [0,T]$.
\end{itemize}
\end{theorem}
The following comparison principle for (\ref{frameeq}) is a straightforward
consequence of Theorem \ref{positivethm} (ii).
\begin{theorem}
Let $v, \tilde{v} \in \Gamma_{T}$ be such that 
$v, \tilde{v}, v_x, \tilde{v}_{x}$ are bounded on 
$\rs \times (0,T]$, $\nop (\tilde{v})(x,t) \leq 0$ and 
$\nop (v)(x,t) \geq 0$ for $(x,t) \in  \rs \times (0,T]$.
Suppose that $\tilde{v}(x,0)-v(x,0) \geq 0$ for all 
$x \in \rs$.
Then $\tilde{v}(x,t) \geq v(x,t)$ for all $(x,t) \in \rs
\times [0,T]$.
\label{compthm}
\end{theorem}
\vspace{4mm}
{\bf Global existence and {\em a priori} bounds}\\[2mm]
The abstract existence theory of \cite[p 253-275]{Lunardi}
applies to the concrete problem
\begin{eqnarray}
v_t & = & A v_{xx} +cv_{x}+ f(v, v_x), ~~ x \in \rs,~ t>0,~ v(x,t) \in \RN,
\label{appeq} \\
v( \cdot , 0) & = & \phi,
\label{appdata}
\end{eqnarray}
where  $A$ satisfies {\bf (a)}, $c \in \rs$ and $f \in \fcont$.

The local existence of a unique  solution of (\ref{appeq}), (\ref{appdata})
and continuous dependence  on the initial data (\ref{appdata})
are a consequence of \cite[p 258, Theorem 7.1.2, p266, Proposition 7.1.9 
and p268, Propostion 7.1.10 and p270, Remark 7.1.12]{Lunardi}.
$\buc^{1}$ is a suitable choice of space between $\buc^{2}$ and $\buc$ for the
initial data $\phi$ -- see \cite[ p 253]{Lunardi},
the embeddings (\ref{halfone}) and (\ref{halfbeta}) and  Lemma \ref{lempartsectorial}. 
The result is the following.

\begin{proposition}
Let $f \in \fcont$ and $\phi \in \cont$. 
Then there exists a maximal $\tau (\phi) \in (0,\infty ]$ such that
there exists a function $V^{\phi} \in C^{1}((0, \tau(\phi)), \buc)
\cap C((0, \tau (\phi)), \buc^{2}) \cap C([0, \tau (\phi)), \buc^{1})$
such that $v^{\phi}$ defined by $v^{\phi}(x,t) = V^{\phi}(t)(x)$
for each $x \in \rs$, $t \in [0, \tau (\phi))$ satisfies
(\ref{appeq}), (\ref{appdata}). Moreover, there is a unique function
$V^{\phi} : [0, \tau(\phi)) \rightarrow \buc^{1}$ with these properties. 

In addition,  given $0<T< \tau(\phi)$, there exist $r, K >0$, depending on
$\phi$ and $T$, such that if $\tilde{\phi} \in \cont$ is such that
$\| \phi - \tilde{\phi} \|_{\cont} < r$,
then $\tau (\tilde{\phi}) \geq T$ and $$\|v^{\phi} \fixt - 
v^{\tilde{\phi}}
\fixt \|_{\cont} \leq K \|\phi - \tilde{\phi} \|_{\cont}
~~\mbox{for each}~~ 0 \leq t \leq T.$$

\label{localthm}
\end{proposition}

Under a growth hypothesis on $f$, the following global existence result, 
 conditional on an 
{\em a priori} bound on $\| v( \cdot, t) \|_{\buc}$, is a consequence of 
\cite[p 266, Proposition 7.1.8, p 268, Proposition 7.1.10 and p 272
Proposition 7.2.2]{Lunardi}.

\begin{proposition}
Suppose that $f \in \fcont$ satisfies the growth condition {\bf (f4)}.
Let $\phi \in \cont$ be such that 
\begin{equation}
\sup_{0 \leq \tilde{t} < \tau(\phi)} \| v^{\phi} (\cdot, \tilde{t} )\|_{\buc}
 = K < 
\infty,
\label{condition}
\end{equation}
where $v^{\phi}$ and $\tau(\phi)$ are as in Proposition \ref{localthm}.
Then $\tau (\phi) = \infty$.
\label{condthm}
\end{proposition}




\noindent {\bf Sub- and supersolutions}\\[2mm]
Theorem \ref{compthm} enables verification of condition (\ref{condition})
under additional hypotheses on $f$ and  $\phi$. Suppose
 that $f \in \fcont$ satisfies {\bf (f1) - (f4)}.
Let $e_{0} = \min_{1 \leq i \leq N} \{ \p_{i} - \m_{i} \} >0$.
 Conditions
{\bf (f2)-(f3)} and the Perron-Frobenius Theorem together imply the existence
of $\nup, \num >0$ and vectors $\eep, \eem \in \RN, \eg^{\pm} >0,
 \| \eg^{\pm} \| =1$ such that $d_{q}f [\epm, 0] \eg^{\pm} =
 - \nu^{\pm} \eg^{\pm}$.
Since $f \in \fcont$, it follows that there exist $ \pz,  \nu >0$, $
\qz \in (0, \half e_{0}), \etaz \in (0, \half e_{0})$
such that
\begin{equation}
\label{fprop}
\left.
\begin{array}{c}
q \in \RN, \|q \| \leq \qz \\ p \in \RN, \| p \| \leq \pz \\ \eta \in 
(0, \etaz) 
\end{array}
\right\}
\Rightarrow 
\left\{ \begin{array}{c}
d_{q}f [ \epm + q -  \eta \eepm, p] \eg^{\pm} < - \nu \eepm, \\ 
d_{q}f [ \epm + q +  \eta \eepm, p] \eg^{\pm} < - \nu \eepm, \\
f(\epm + q -  \eta \eepm , p) - f( \epm + q, p) \geq  \nu \eta \eepm,\\
f(\epm + q +  \eta \eepm , p) - f( \epm + q, p) \leq  \nu \eta \eepm.
\end{array} \right. 
\end{equation}
Suppose that \tw ~holds. 
 Since $w'(x) \rightarrow 0$ as $|x| \rightarrow \infty$,
we can choose $\qz, \etaz$ smaller if necessary to ensure that
\begin{equation}
\label{equilibrium}
\| w(x) - \epm \| < \qz + \etaz \Rightarrow \| w'(x) \| < \pz.
\end{equation}
Let  $\pg \in C^{\infty}(\RN, \RN)$ be such that ${\mathfrak{p}}_{i}(q) =
\pgi (q_{i})$ for each $i \in \{ 1, \ldots,
N \}$, $q \in \RN$ (the $i$-th component of $\pg$ depends only on the $i$-th
component of its argument), where  
 $\pgi \in C^{\infty}(\rs, \rs)$
is a smooth monotone function with
\begin{equation}
\label{pg}
\pgi (\omega) = \eep_{i} ~\mbox{when}~ | \p_{i} - \omega | \leq \qz, 
~~\mbox{and}~~ \pgi (\omega) = \eem_{i} ~\mbox{when}~ | \m_{i} - \omega | \leq
\qz.
\end{equation}

The following construction of sub- and super-solutions is an extension,
to nonlinearities $f$ depending on $v_{x}$, of constructions in \cite{FM}
and \cite{Hutson}.

\begin{theorem}
\label{fronttrapthm}
There exist $\alpha_{0} >0$ and $\hat{\eta} \in (0, \etaz]$ such that for any
$ x_{0}, x_{1} \in \rs$ and any $\eta \in [0, \hat{\eta}]$,
$$ \nop (\sub) (x,t) \geq 0 ~~\mbox{and}~~ \nop (\super) (x,t) \leq 0
 ~~\mbox{for all}~~ x \in \rs, t \geq 0, $$
where for each $i \in \{ 1, \ldots , N \}$,
\begin{equation}
\label{subdef}
(\sub)_{i}(x,t) = w_{i}(x-x_{0} + \eta \alpha_{0} e^{- \nu t})
 - \eta e^{- \nu t}\pgi ( w_{i}(x-x_{0} + \eta \alpha_{0} e^{- \nu t}))
\end{equation}
and
\begin{equation}
\label{superdef}
(\super)_{i}(x,t) = w_{i}(x + x_{1} - \eta \alpha_{0} e^{ -\nu t})
+ \eta e^{ -\nu t} \pgi (w_{i}(x + x_{1} - \eta \alpha_{0} e^{ -\nu t}) ).
\end{equation}
Here the $c$ in (\ref{defN}) 
 is the velocity of the
wave $w$.
\end{theorem}

\proof Let $x_{0}, x_{1} \in \rs$ be arbitrary.
Let $\alpha_{0} >0$ (to be fixed later), and let $\eta \in (0, \etaz]$.
  Define $\sub$ and $\super$ as in
(\ref{subdef}) and (\ref{superdef}). We will prove the result for $\sub$; the
argument for $\super$ is similar. 

Unless otherwise indicated, $w, w'$ are to be evaluated at the point $\argxz$.
Fix $t \geq 0$. First let $x$ be such that $\| \p - w(x-x_{0}+\eta \alpha_{0}
e^{- \nu t}) \| \leq \qz /2$. For such $x$, 
$\pgi'( w_{i} \argxz ) =0$ and $\pgi (w_{i} \argxz) =\eep_{i}$ for each 
$i \in \{ 1, \ldots, N\}$. Hence (\ref{fprop}), (\ref{equilibrium})
 together with the facts that
$w$ is a stationary solution of (\ref{frameeq}) and that
$w'(s)>0$ for all $s$ yield that 
\begin{eqnarray*}
\nop (\sub) (x,t) & = & \nu \eta \alpha_{0} \ent w'   -
\nu \eta \ent \eep   
+ f(w   - \eta \ent \eep, w'  ) 
 - f(w  , w'  ) \\  
& \geq & \nu \eta \ent \eep - \nu \eta \ent \eep = 0.
\end{eqnarray*}
Similarly, $\nop (\sub )(x,t) \geq 0$ when $\| \m - w \argxz  \| \leq
\qz/2$.

Now let $x \in \rs$ be such that $$\| \m - w \argxz  \| \geq \qz/2
~~\mbox{and}~~
\| \p - w \argxz \| \geq \qz/2$$.
Since $w' >0$, there exists $\beta >0$, depending only on $w$ and $\qz$,
such that for each $i \in \{ 1, \ldots , N \} $,
\begin{equation}
\label{lowerderivbound}
\| \m - w(s) \| \geq \frac{\qz}{2} ~\mbox{and}~ 
\| \p - w(s) \| \geq \frac{\qz}{2}  \Rightarrow w_{i}'(s) \geq \beta.
\end{equation}
Let $i \in \{ 1, \ldots, N\}$.
Since $w$ is a stationary solution of (\ref{frameeq}),
$$
\nop_{i} (\sub) \xt  =  \nu \eta \alpha_{0} \ent w_{i}'  
- \nu \eta \ent \pgi(w_{i}  ) 
 - \nu \eta^{2} \alpha_{0} e^{-2 \nu t}
\pgi' (w_{i})   w_{i}' 
 - \eta \ent A_{i} [ \pgi'' (w_{i}) (w_{i}')^{2} 
+ \pgi'(w_{i})  w_{i}'' ]  $$ $$   - \eta \ent c \pgi'(w_{i}) w_{i}' 
 + f_{i}(w - \pg \eta \ent, w' - \eta \ent d \pg [w] w') 
- f_{i}( w, w'). 
$$
By the Mean Value Theorem and the properties  of $\pg$ and $w$,
\begin{equation}
\label{boundedvector}
f_{i}(w - \pg \eta \ent, w' - \eta \ent d \pg [w] w') 
- f_{i}( w, w')
= {\mathfrak{q}}_{i} \xt \eta \ent,
\end{equation}
where ${\mathfrak{q}}_{i} \xt$ is bounded independently of $x \in \rs$ and
$t \geq 0$. So
\begin{eqnarray}
\nop_{i} (\sub) \xt & = & \eta \ent \left\{ {\mathfrak{q}}_{i} \xt - 
\nu \pgi (w_{i}) - A_{i} [\pgi''(w_{i}) (w_{i}')^{2} + \pgi'(w_{i}) w_{i}'']
- c \pgi'(w_{i}) w_{i}' \right\} \nonumber \\
& & \hspace*{2cm}
\mbox{} + \nu \eta \alpha_{0} \ent w_{i}' \{ 1 - \eta \ent \pgi'(w_{i}) \}.
\label{Ndominate}
\end{eqnarray}
Since $d \pg [\cdot]$ is uniformly bounded, there exists $\hat{\eta} \in (0, \etaz]$
such that
\begin{equation}
\label{etahat}
\eta \in (0, \hat{\eta}] \Rightarrow 1 - \eta \ent | \pgi'(\omega) | \geq \half
~~\mbox{for each}~~ \omega \in \rs.
\end{equation}
(We need that $1 + \eta \ent \pgi' (\omega) \geq \half$ for the analysis of
 $\super$.)
So since $w_{i}'$ satisfies (\ref{lowerderivbound}),
\begin{eqnarray*}
\lefteqn{\nop_{i} (\sub) \xt  \geq } \\
& &  \eta \ent \left\{ {\mathfrak{q}}_{i} \xt - 
\nu \pgi (w_{i}) - A_{i} [\pgi''(w_{i}) (w_{i}')^{2} + \pgi'(w_{i}) w_{i}''] 
- c \pgi'(w_{i}) w_{i}' + \half \nu \beta \alpha_{0} \right\}.
\end{eqnarray*}
Whence we can choose $\alpha_{0} >0$, dependent  on $\pg$ and 
$w$ but independent of $x$ and $t$,
 such that $\nop_{i} (\sub) \xt \geq 0$. The result follows. \qed

\begin{theorem}
\label{frontglobalthm}
Suppose that $f \in \fcont$ satisfies {\bf (f1) - (f4)}. Then there
exists $\hat{\eta} >0$ such that if $\phi \in \buc^{1}$ is such that
there exists $ \eta \in (0, \hat{\eta})$ with
\begin{equation}
\label{trap}
\m - \eta \eem \leq \phi (x) \leq \p + \eta \eep ~~\mbox{for all}~~ x \in \rs,
\end{equation}
and
\begin{equation}
\label{limsup}
\limsup_{x \rightarrow \infty} | \phi_{i} (x) - \p_{i} | \leq \eta \eep_{i},~~
 \limsup_{x \rightarrow - \infty} | \phi_{i} (x) - \m_{i}|
\leq \eta \eem_{i} ~~\mbox{for each}~~ i \in \{ 1, \ldots, N\},
\end{equation}
then $\tau (\phi) = \infty$, and there exist $\xzp, \xop \in \rs$ such that
\begin{equation}
\label{timetrap}
{\mathbf{s}}_{\hat{\eta}, \xzp} \xt \leq v^{\phi} \xt \leq
 {\mathbf{S}}_{\hat{\eta}, \xop} \xt ~~\mbox{for all}~~ x \in \rs, t \geq 0.
\end{equation}
\end{theorem}

\proof Let $\hat{\eta}$ be as in Theorem \ref{fronttrapthm}, and let
$\phi \in \buc^{1}$ satisfy (\ref{trap}, \ref{limsup}) for some $\eta \in 
(0, \hat{\eta})$. Now given $x_{0}, x_{1} \in \rs$, 
\begin{eqnarray}
{\mathbf{s}}_{\hat{\eta}, x_{0}} (x,0) &= & w (x-x_{0} + \hat{\eta}
 \alpha_{0}) - \hat{\eta} \pg (w(x-x_{0} + \hat{\eta} \alpha_{0}))
\label{subzero} \\
{\mathbf{S}}_{\hat{\eta}, x_{1}} (x,0) & = & w(x+ x_{1} - \hat{\eta}
\alpha_{0})+ \hat{\eta} \pg (w(x+x_{1} - \hat{\eta} \alpha_{0}))
\label{superzero}
\end{eqnarray}
for each $x \in \rs$. Recall (\ref{pg}).
 So from (\ref{trap}), (\ref{limsup}),
(\ref{subzero}), (\ref{superzero}) and the fact that $\eta < \hat{\eta}$,
it follows
that there exist $\xzp, \xop \in \rs$ such that
\begin{equation}
{\mathbf{s}}_{\hat{\eta}, \xzp} (x,0) \leq \phi (x) \leq 
{\mathbf{S}}_{\hat{\eta}, \xop} (x,0) ~~\mbox{for all}~~ x \in \rs.
\end{equation}
This, together with Theorem \ref{fronttrapthm}, allow application of
Theorem \ref{compthm} to get that
\begin{equation}
\label{vphitrap}
{\mathbf{s}}_{\hat{\eta}, \xzp} \xt \leq v^{\phi} \xt \leq
{\mathbf{S}}_{\hat{\eta}, \xop} \xt ~~\mbox{for all}~~ x \in \rs, 
0 \leq t < \tau( \phi).
\end{equation}
Whence condition (\ref{condition}) is satisfied. The result follows from 
Proposition \ref{condthm}. \qed 

The wave-dependent sub- and super-solutions constructed above are valuable
in analysing the stability of the wave $w$. The following is another,
simple but useful, route to verification of condition (\ref{condition}).

\begin{theorem}
\label{globalthm}
Suppose that $f \in \fcont$ satisfies {\bf (f1)-(f4)}. Then there exists
$\hat{\eta} >0$ such that if $\phi \in \buc^{1}$ is such that there exists
$\eta \in [0, \hat{\eta}]$ such that
\begin{equation}
\label{databound}
\m - \eta \eem \leq \phi (x) \leq \p + \eta \eep ~~\mbox{for all}~~ x \in \rs,
\end{equation}
then $\tau(\phi) = \infty$, and
\begin{equation}
\label{solnbound}
\m - \eta \eem \leq v^{\phi} \xt \leq \p + \eta \eep ~~\mbox{for all}~~
x \in \rs, t \geq 0.
\end{equation}
\end{theorem}

\proof Let $\hat{\eta}$ be as in Theorem \ref{fronttrapthm} and let $\phi \in
\buc^{1}$ satisfy (\ref{databound}) for some $\eta \in [0, \hat{\eta}]$.
Then since $f( \p, 0) = f(\m, 0)=0$ (by {\bf (f3)}), it follows from
(\ref{fprop}) that $f(\m - \eta \eem, 0) >0, f(\p + \eta \eep, 0) <0$.
So with $\nop$ as defined in (\ref{defN}), $\nop(\m - \eta \eem, 0) >0$
and $\nop (\p + \eta \eep, 0) <0$. It then follows from Theorem \ref{compthm}
that 
\begin{equation}
\m - \eta \eem \leq v^{\phi} \xt \leq \p + \eta \eep ~~\mbox{for}~~
0 \leq t \leq \tau(\phi).
\end{equation}
Whence condition (\ref{condition}) is satisfied. The result follows from
Proposition \ref{condthm}. \qed

\noindent {\bf{\em A priori} bounds}\\[2mm]
The  derivatives of $v^{\phi}$ can be estimated independently of
 the exact choice of $\phi$ satisfying (\ref{databound}), as follows.

\begin{theorem}
Let $f, \hat{\eta}$ be as in Theorem \ref{globalthm} and let $t_{0} >0$.
 Then there exists
$K(t_{0}) >0$
 such that if $\phi \in \cont$ satisfies (\ref{databound}) for
some $\eta \in [0, \hat{\eta}]$, then for all $t \geq t_{0}$,
\begin{equation}
\| v^{\phi} \fixt \|_{\buc^{2}} \leq K(t_{0}).
\label{xbound}
\end{equation}
\label{boundthm}
\end{theorem}

\proof Since $f$ satisfies {\bf (f1), (f4)} and (\ref{solnbound}) holds,
 the single-equation analysis of
\cite[Chapter V, \S3, p 437, Theorem 3.1]{LSU} implies the existence
of $K_{1}(t_{0})>0$ such that $\| v^{\phi} \fixt \|_{\buc^{1}}
 \leq K_{1}(t_{0})$ for all $t \geq t_{0}$.
This enables application of \cite[Chapter VII, \S 5, p 586, Theorem 5.1]{LSU}
to obtain (\ref{xbound}).  \qed \vspace{2mm}

\bigskip
\bigskip
\noindent \textsc{E.C.M.Crooks $^*$}\\
Department of Mathematical Sciences\\
University of Bath\\
Bath BA2 7AY \\
United Kingdom
\vspace{2mm} \\
$^*$ Current address : Balliol College, Oxford, OX1 3BJ, U.K.,
\texttt{crooks@maths.ox.ac.uk} 

\end{document}